\documentclass[a4paper,12pt,twoside]{article}
\usepackage[english]{babel}
\usepackage[utf8]{inputenc}

\usepackage{latexsym}
\usepackage{amssymb}
\usepackage{mathtools}
\usepackage{amsthm}
\usepackage{amsfonts}
\usepackage{graphicx}
\usepackage{epstopdf}
\usepackage{color}
\usepackage{enumitem}

\usepackage{aliascnt}
\usepackage[capitalise]{cleveref}

\allowdisplaybreaks

\numberwithin{equation}{section}

\newtheorem{theorem}{Theorem}[section]
\newaliascnt{lemma}{theorem}
\newtheorem{lemma}[lemma]{Lemma}
\aliascntresetthe{lemma}
\newaliascnt{corollary}{theorem}
\newtheorem{corollary}[corollary]{Corollary}
\aliascntresetthe{corollary}
\newaliascnt{proposition}{theorem}

\aliascntresetthe{proposition}
\newaliascnt{example}{theorem}
\newtheorem{example}[example]{Example}
\aliascntresetthe{example}
\newaliascnt{remark}{theorem}

\aliascntresetthe{remark}
\newaliascnt{assumption}{theorem}

\aliascntresetthe{assumption}
\newaliascnt{problem}{theorem}
\newtheorem{problem}[problem]{Problem}
\aliascntresetthe{problem}
\theoremstyle{definition}
\newaliascnt{definition}{theorem}
\newtheorem{definition}[definition]{Definition}
\aliascntresetthe{definition}

\author{Ferenc Weisz
\\
ELTE E\"otv\"os Lor\'and University \\
Department of Numerical Analysis \\
H-1117 Budapest, P\'azm\'any P. s\'et\'any 1/C., Hungary \\
e-mail: weisz@inf.elte.hu
\and
Guangheng Xie
\\
School of Mathematics and Statistics, HNP-LAMA \\
Hunan Research Center of the Basic Discipline\\
for Analytical Mathematics, Central South University \\
Changsha 410083, the People's Republic of China \\
e-mail: xieguangheng@csu.edu.cn
\and
Dachun Yang (corresponding author)
\\
Laboratory of Mathematics and Complex Systems\\
(Ministry of Education of China) \\
School of Mathematical Sciences, Beijing Normal University \\
Beijing 100875, the People's Republic of China \\
e-mail: dcyang@bnu.edu.cn}

\title{A characterization of almost everywhere convergence of
subsequences of Walsh--Fourier partial sums}

\date{}

\begin{document}
\maketitle

\begin{abstract}
We characterize the increasing sequences $(\nu_j)$ in the set of positive
integers $\mathbb N$ along which the Walsh--Fourier partial sums
$s_{\nu_j}(f)$ of every integrable function $f$ converge almost everywhere
to $f$. This resolves
a problem that Konyagin posed in his invited address at the 2006
International Congress of Mathematicians~\cite{K06}. The criterion concerns the
binary expansions of the integers $\nu_j$. Divide each $\nu_j$ by $2^k$.
Convergence holds if and only if, for some $K\in\mathbb N$, the number of
distinct positive quotients whose remainder has binary variation larger
than $K$ is bounded independently of the scale $k$. This is also
equivalent to the weak type $(1,1)$ inequality for
$\sup_{j\in\mathbb N}|s_{\nu_j}(f)|$, and if the criterion fails, then the
partial sums of some integrable function along $(\nu_j)$ are unbounded at
every point. The criterion contains the
bounded variation condition and Konyagin's condition for separated
spectra, and we give examples, for convergence and for divergence, that
are not covered by the earlier results.

\vspace{5mm}

\vspace{3mm}
\begin{sloppypar}
\end{sloppypar}
\end{abstract}

\section{Introduction}\label{s1}
The Walsh--Paley system is the dyadic analogue of the trigonometric
system \cite{Walsh,Paley,Fine,SWS}. Billard \cite{Billard} proved that
the Walsh--Fourier partial sums of every $f\in L_2[0,1)$ converge to
$f$ almost everywhere. Sj\"olin \cite{Sjolin} extended this result to
$f\in L_p[0,1)$ for $p\in(1,\infty)$. These are the Walsh counterparts of
the theorems of Carleson \cite{Carleson} and Hunt \cite{Hunt} for
trigonometric Fourier series. Thiele \cite{Thiele00} gave a proof
using time--frequency analysis of Walsh operators; see also
\cite{Demeter,LaceyThiele}.
At the $L_1$ endpoint, Stein \cite{Stein} and Schipp
\cite{Schipp67,Schipp69} gave Walsh analogues of Kolmogorov's
divergence examples. Their examples show that the Walsh--Fourier
partial sums of an integrable function can diverge almost everywhere,
and even everywhere.

On the other hand, the
partial sums with index $2^k$ are the conditional expectations of the
function with respect to the dyadic $\sigma$-algebras, so that they
converge almost everywhere for every integrable function by the
martingale convergence theorem.

This leads to a question that has no counterpart for trigonometric
series: along which fixed sequences of indices do the partial sums of
every integrable function converge almost everywhere? For trigonometric
series there is no such sequence. By a theorem of Gosselin
\cite{Gosselin}, every increasing sequence admits an integrable function
whose partial sums along it are unbounded almost everywhere. Totik
\cite{Totik} showed that the unboundedness can be achieved at every
point, and Konyagin \cite{K05} carried Totik's theorem over to the
Orlicz classes $\varphi(L)$ with $\varphi(u)=o(u\log\log u)$. For Walsh
series the sequence $(2^k)_{k\in\mathbb N}$ shows that such sequences exist, and
Konyagin raised the following problem in his invited address at the
International Congress of Mathematicians in 2006; see also Do and Lacey
\cite[p.~253]{DL}.

\begin{problem}[Konyagin {\cite[Problem~3.3, p.~1398]{K06}}]\label{pr1.1}
Find a necessary and sufficient condition on a sequence $(\nu_j)_{j\in\mathbb N}$ in
$\mathbb N$ under which the Walsh--Fourier partial sums
$s_{\nu_j}(f)$ converge to $f$ almost everywhere for every integrable function
$f$ on $[0,1)$.
\end{problem}

Problem \ref{pr1.1} is also stated as open in the very recent paper \cite[Remark~1]{GM26} 
by Goginava and Mukhamedov. Three classes of sequences have been
understood so far. The classical sufficient condition is bounded
variation. If the binary variation is uniformly bounded along the sequence, then the partial
sums along it converge almost everywhere for every integrable function
\cite{Fine,SWS}. Konyagin \cite{K93} answered a question of Balashov by showing that bounded variation 
is not a necessary condition. His result guarantees almost everywhere convergence of the partial sums 
for every integrable function along any sequence with separated spectra, 
and sequences of this type may have unbounded variation. In the
opposite direction, Goginava and Oniani \cite{GO} proved that if the
spectra are nested and the variation is unbounded, then there is an integrable function whose
partial sums along the sequence are unbounded everywhere. Their
construction goes back to Schipp \cite{Schipp67,Schipp69} and Konyagin
\cite{K05}. For the square partial sums of double Walsh--Fourier series, the situation is
different. There, bounded variation of the sequence is necessary for
almost everywhere convergence for every integrable function \cite{K93},
but it is not sufficient \cite{Luk95}. Hence only sequences with bounded
variation have to be examined, and among them Lukomskij found a necessary
and sufficient condition for double series \cite{Luk95}, and extended it
to every dimension greater than one \cite{Luk98}. Convergence along 
lacunary sequences under integrability assumptions 
slightly stronger than \(L_1\) has been investigated by Do and Lacey \cite{DL}, 
Lie \cite{Lie12,Lie13,Lie,Lie19}, and Di Plinio \cite{DP}. For the arithmetic averages,
$N^{-1}\sum_{j=1}^{N}s_{\nu_j}(f)$,
G\'at proved almost everywhere convergence for every integrable function
$f$ along every lacunary sequence \cite{Gat10}.
This result was subsequently extended to sequences satisfying a 
weaker growth condition in \cite{Gat19,Gog26}. On the negative side, Blahota and G\'at \cite{BG26} 
established that such averages can diverge almost everywhere 
for a suitable sequence and a suitable integrable function.

In this paper we solve Konyagin's problem, that is, Problem \ref{pr1.1}. Our criterion looks at the
binary expansions of the terms of the sequence at every dyadic scale.
Dividing a term by $2^k$, we obtain a quotient and a remainder, and the
remainder consists of the lowest $k$ binary digits of the term. We prove
that the partial sums along the sequence converge almost everywhere for
every integrable function if and only if there exists $K\in\mathbb N$ such
that, at every scale, only boundedly many distinct positive quotients
occur together with a remainder of binary variation larger than $K$.
The bound may depend on the sequence, but not on the scale. The precise
statement is \cref{t3.2}, and \cref{c3.3} gives an equivalent form with a
single parameter. The same condition characterizes the sequences for
which the maximal operator of the selected partial sums is of weak type
$(1,1)$. The classical conditions are special cases: for bounded
variation no quotient is counted at all, and for separated spectra at
most one quotient is counted at each scale.

For Vilenkin systems, Areshidze, Persson and Tephnadze
\cite[Theorem~1]{APT25} give a sufficient condition for almost everywhere
convergence in terms of variation. They ask for a characterization when
the variation is unbounded \cite[Open Problem~3]{APT25}. We answer this
question in the Walsh--Paley case.

If the criterion fails, then the divergence is as strong as it can be:
there is an integrable function whose partial sums along the sequence
are unbounded at every point (\cref{t3.6}). Consequently, for every
sequence exactly one of two alternatives occurs. Either the partial sums
along the sequence converge almost everywhere for every integrable
function, or they diverge everywhere for some integrable function. This
extends the theorem of Goginava and Oniani \cite{GO} on nested spectra to
arbitrary sequences. For the full sequence of partial sums and for
homogeneous Banach spaces of integrable functions, alternatives of this
kind are due to Katznelson \cite{Katz66} in the trigonometric case and to
Harris and Wade \cite{HW78} in the Walsh case, and Oniani \cite{On22}
extended them to systems of characters of compact abelian groups.

In \cref{s7}, we recover the theorems on bounded variation, on separated spectra and on
nested spectra. We then construct a sequence along which the partial sums of every
integrable function converge almost everywhere, although it cannot be
covered by finitely many sequences of bounded variation or with separated
spectra. We also give a sequence along which the criterion fails, so that
the partial sums of some integrable function are unbounded at every
point, although none of its subsequences is at bounded distance from a sequence with
nested spectra and the results of \cite{GO} do not apply.

The remainder of the paper is organized as follows. In \cref{s2},
we collect the dyadic identities used in the proofs and recall the
maximal principle. \Cref{s3} contains the main results and the
deduction of the convergence criterion from the two estimates in
\cref{t3.4}. In \cref{s4}, we prove the upper bound for the weak type
$(1,1)$ constant using a Walsh version of the multifrequency
Calder\'on--Zygmund decomposition of Nazarov, Oberlin and Thiele
\cite[Theorem~1.1]{NOT}. In \cref{s5}, we prove the lower bound using
a product of Walsh polynomials. In \cref{s6}, we prove \cref{c3.5}
and \cref{t3.6}, which give unbounded partial sums almost everywhere
and everywhere, respectively. Finally, in \cref{s7}, we recover
earlier results on bounded variation, separated spectra, and nested
spectra, and give two examples illustrating convergence and divergence.

We conclude the introduction with some notational conventions.
Let $\mathbb N\coloneqq\{1,2,\ldots\}$ and
$\mathbb Z_+\coloneqq\{0,1,2,\ldots\}$ denote the sets of
\emph{positive integers} and \emph{nonnegative integers}, respectively.
For any set $A$, let $\#A$ denote its cardinality, with $\#A=\infty$
if $A$ is infinite. For any real number $u$, let $\lfloor u\rfloor$ and
$\lceil u\rceil$ denote the greatest integer at most $u$ and the least
integer at least $u$, respectively. We write $X\lesssim Y$ if $X\leq CY$
for an absolute constant $C$. Empty sums are zero and empty products
are one. We work with real-valued functions. The convergence results
extend to complex-valued functions by applying them to the real and
imaginary parts separately.

\section{Preliminaries}\label{s2}
In this section, we prepare for the proof of the convergence criterion
in our main result. To this end, we collect the dyadic identities for
Walsh--Fourier partial sums and use Stein's maximal principle
\cite{Stein} to reduce almost everywhere convergence for every
integrable function to a weak type $(1,1)$ estimate for the associated
maximal operator.

\subsection{The Walsh--Paley system}
We work with the Walsh--Paley system on $[0,1)$ \cite{Walsh,Paley}; see, for example, \cite{Fine,SWS,wk2}. 
For any Lebesgue
measurable set $A\subset[0,1)$, let $\lambda(A)$ denote its
\emph{Lebesgue measure}. For any interval $I$, let $|I|$ denote its
\emph{length}. For any $x\in[0,1)$, write its \emph{binary expansion} as
\[
 x=\sum_{l=0}^{\infty}x_l2^{-l-1}\qquad(x_l\in\{0,1\}).
\]
If $x$ is a dyadic rational, we choose the expansion that is eventually zero.

Let $G\coloneqq\{0,1\}^{\mathbb Z_+}$ be the \emph{dyadic group}, a compact
group with coordinatewise addition modulo two and the product Haar probability
measure. Let $G_0\subset G$ consist of the sequences that are not
eventually equal to one. The complement $G\setminus G_0$ is countable.
Define the \emph{binary digit map} $\beta\colon[0,1)\to G_0$ by
$\beta(x)\coloneqq(x_l)_{l\in\mathbb Z_+}$. The convention on binary
expansions makes $\beta$ a bijection. The map preserves measure because
the inverse image of a cylinder specified by the first $k$ digits is a
dyadic interval of the same measure $2^{-k}$. For any $x,y\in[0,1)$,
we define their \emph{dyadic sum} by
\[
 x\dot+y\coloneqq
 \begin{cases}
  \beta^{-1}\bigl(\beta(x)+\beta(y)\bigr),
     &\beta(x)+\beta(y)\in G_0,\\
  0,&\beta(x)+\beta(y)\notin G_0.
 \end{cases}
\]
For fixed $y\in[0,1)$, the second case occurs for only countably many
$x$. Since translation by $\beta(y)$ preserves Haar measure, the map
$x\mapsto x\dot+y$ preserves Lebesgue measure. In the first case,
$(x\dot+y)\dot+y=x$ and the binary digits of $x\dot+y$ are $|x_l-y_l|$.
If $y$ is dyadic, then $\beta(y)$ has only finitely many nonzero digits,
so the first case applies for every $x$.

Every $n\in\mathbb Z_+$ has the binary expansion
$
 n=\sum_{l=0}^{\infty}n_l2^l$ ($n_l\in\{0,1\}$).
For any $m,n\in\mathbb Z_+$, we write $m\oplus n$ for the integer
with the binary digits $|m_l-n_l|$. The \emph{Rademacher functions} and the
\emph{Walsh--Paley functions} are
\[
 r_l(x)\coloneqq (-1)^{x_l}\quad(l\in\mathbb Z_+,\ x\in[0,1)),\qquad
 w_n\coloneqq \prod_{l=0}^{\infty}r_l^{n_l}\quad(n\in\mathbb Z_+).
\]
Only finitely many factors differ from one, and $w_0=1$. The Walsh system
$(w_n)_{n\in\mathbb Z_+}$ is orthonormal, and $w_mw_n=w_{m\oplus n}$ for
$m,n\in\mathbb Z_+$. Moreover, if $y\in[0,1)$ is fixed, then
$w_n(x\dot+y)=w_n(x)w_n(y)$ for almost every $x\in[0,1)$. 
The finite linear combinations of the functions $w_m$ are called
\emph{Walsh polynomials}. These polynomials are dense in $L_1[0,1)$
\cite[p.~12]{SWS}. For any $f\in L_1[0,1)$, its \emph{Walsh--Fourier
coefficients} and \emph{partial sums} are
\begin{align*}
 \widehat f(m)\coloneqq \int_0^1 f w_m\,d\lambda\quad(m\in\mathbb Z_+),\qquad
 s_n(f)\coloneqq \sum_{m=0}^{n-1}\widehat f(m)w_m\quad(n\in\mathbb N),\qquad
 s_0(f)\coloneqq 0.
\end{align*}
For any $p\in[1,\infty),$ we write
$\|f\|_{L_p[0,1)}\coloneqq (\int_0^1|f|^p\,d\lambda)^{1/p}$, and
$\|f\|_{L_\infty[0,1)}$ denotes the \emph{essential supremum} of $|f|$.

The \emph{binary variation} of $n\in\mathbb Z_+$ and the set of its \emph{occupied
binary positions} are
\begin{align*}
 V(n)\coloneqq n_0+\sum_{l=1}^{\infty}|n_l-n_{l-1}|,\qquad
 \operatorname{Sp}(n)\coloneqq \{l\in\mathbb Z_+:\ n_l=1\}.
\end{align*}
Thus $V(0)=0$, and for $n\in\mathbb N$ the number $V(n)$ counts the changes
between consecutive terms of the sequence $0,n_0,n_1,\ldots$.
Since this sequence begins with $0$ and ends with zeros, $V(n)$ is even;
for example, $V(2^l)=2$ for $l\in\mathbb Z_+$. The \emph{remainder} of
$n$ upon division by $2^k$ is
\[
 n\bmod2^k\coloneqq n-2^k\lfloor n/2^k\rfloor=\sum_{l=0}^{k-1}n_l2^l .
\]

\subsection{Averages on dyadic intervals}
A \emph{dyadic interval} is an
interval of the form $[i2^{-k},(i+1)2^{-k})$ with $k\in\mathbb Z_+$ and
$i\in\{0,1,\ldots,2^k-1\}$. For $k\in\mathbb Z_+$ and $x\in[0,1)$, the dyadic interval
of length $2^{-k}$ that contains $x$ is
\[
 I_k(x)\coloneqq \{y\in[0,1):\ y_l=x_l\ (0\leq l<k)\},\qquad |I_k(x)|=2^{-k}.
\]
In particular $I_0(x)=[0,1)$, and we put $I_k\coloneqq I_k(0)=[0,2^{-k})$. 
For an integrable function $h$, the \emph{conditional expectation} with respect
to the $\sigma$-algebra generated by the dyadic intervals of length
$2^{-k}$ is
\begin{align*}
 E_k(h)(x)\coloneqq 2^k\int_{I_k(x)}h\,d\lambda\qquad(x\in[0,1)).
\end{align*}
The sequence $(E_k(h))_{k\in\mathbb Z_+}$ is a martingale. We write
\[
 d_k(h)\coloneqq E_k(h)-E_{k-1}(h)\quad(k\in\mathbb N),\qquad
 M(h)\coloneqq \sup_{k\in\mathbb Z_+}E_k(|h|)
\]
for its \emph{martingale differences} and for the \emph{dyadic maximal function} of
$h$. For an interval $I,$ we also use the \emph{average}
\[
 \langle h\rangle_I\coloneqq \frac{1}{|I|}\int_I h\,d\lambda.
\]
We need the next two well-known results from martingale theory (see \cite{wk}):
\begin{align}
 \lambda\{M(h)>\varrho\}&\leq\frac{\|h\|_{L_1[0,1)}}{\varrho}\quad(\varrho>0), &
 \lim_{k\to\infty}E_k(h)&=h\quad\text{almost everywhere}.\label{e2.1}
\end{align}

\subsection{Walsh products, kernels, and binary coefficients}

The \emph{Dirichlet kernels} are $D_n\coloneqq \sum_{m=0}^{n-1}w_m$ $(n\in\mathbb N)$ and
$D_0\coloneqq 0$. For each $x\in[0,1)$, the identity
$w_m(x\dot+t)=w_m(x)w_m(t)$ holds for almost every $t$. Summing over
$m<n$ and using $|D_n|\leq n$, we find that, for any $h\in L_1[0,1)$,
\begin{align*}
 s_n(h)(x)=\int_0^1 h(t)D_n(x\dot+t)\,dt\quad(x\in[0,1)),\qquad
 \|s_n(h)\|_{L_\infty[0,1)}\leq n\|h\|_{L_1[0,1)}.
\end{align*}
For any $y\in[0,1)$, by the change of variables $t\mapsto t\dot+y$, we have
\[
 \widehat{h(\cdot\dot+y)}(m)
 =\int_0^1 h(t)w_m(t\dot+y)\,dt
 =w_m(y)\widehat h(m).
\]
See also \cite[Theorem~II, p.~380]{Fine}.
Inserting this identity into the finite Walsh sum yields
\begin{align}\label{e2.2}
 s_n\bigl(h(\cdot\dot+y)\bigr)(x)=s_n(h)(x\dot+y)
 \quad\text{for almost every }x\in[0,1).
\end{align}

It is well known \cite{SWS,wk2} that
\begin{align}\label{e2.3}
 D_{2^k}(x)=\begin{cases}2^k,&x\in[0,2^{-k}),\\0,&x\in[2^{-k},1)\end{cases}
 \qquad(k\in\mathbb Z_+),
\end{align}
and that, for any $k\in\mathbb Z_+$, $s_{2^k}(h)=E_k(h)$ for any $h\in L_1[0,1)$. For
$q,a\in\mathbb Z_+$ with $a<2^k$,
\begin{align}\label{e2.4}
 D_{2^kq+a}=D_{2^kq}+w_{2^kq}D_a,\qquad
 D_{2^kq}=D_{2^k}\sum_{j=0}^{q-1}w_{2^kj}.
\end{align}
In particular, by \eqref{e2.3},
\begin{align}\label{e2.5}
 D_{2^kq}(z)=0\quad\text{if }z\notin I_k.
\end{align}
For any $k\in\mathbb Z_+$ and
$x\in[0,1)$, the number $2^kx-\lfloor2^kx\rfloor$ has the binary digits
$x_k,x_{k+1},\ldots$, which gives
\begin{align}\label{e2.6}
 w_{2^kq}(x)=w_q\bigl(2^kx-\lfloor2^kx\rfloor\bigr)\qquad(q\in\mathbb Z_+).
\end{align}
For any $q,q'\in\mathbb Z_+$, we have $w_{2^kq}w_{2^kq'}=w_{2^k(q\oplus q')}$. If
$I=[i2^{-k},(i+1)2^{-k})$ is a dyadic interval of length $2^{-k}$, then
$2^kx-\lfloor2^kx\rfloor=2^kx-i$ for $x\in I$, and the substitution
$u=2^kx-i$ maps $I$ onto $[0,1)$. Hence, by \eqref{e2.6} and the orthonormality of the Walsh system,
\begin{align}\label{e2.7}
 \frac1{|I|}\int_Iw_{2^kq}w_{2^kq'}\,d\lambda
    =\int_0^1w_{q\oplus q'}\,d\lambda
    =\begin{cases}1,&q=q',\\0,&q\ne q'\end{cases}
 \qquad(q,q'\in\mathbb Z_+).
\end{align}
In particular $E_k(w_{2^kq})=0$ for every $q\in\mathbb N$.

The following two kernel expansions will be used separately:
\begin{align}\label{e2.8}
 w_nD_n=\sum_{l=0}^{\infty}n_lr_lD_{2^l}
       =\sum_{l=0}^{\infty}c_l(n)D_{2^l},
\end{align}
where the integers $c_l(n)$ are
$
 c_0(n)\coloneqq -n_0$ and
 $c_l(n)\coloneqq n_{l-1}-n_l$ for every $l\in\mathbb N$.
Both expansions are known. The first one is
\cite[Theorem~8\,(iii), p.~28]{SWS}, see also \cite[Lemma~1.6.5]{wk6}.
The second expansion is obtained from the first by Abel's transformation,
since $r_lD_{2^l}=D_{2^{l+1}}-D_{2^l}$ by \eqref{e2.4} with $q=2$, see
\cite[(58), p.~35]{SWS}.

For any $n\in\mathbb Z_+$ and $x\in [0,1),$ applying \eqref{e2.8} at $x\dot+t$, using $w_n(x\dot+t)=w_n(x)w_n(t)$ and $s_{2^l}=E_l$, we deduce that
\begin{align*}
 s_n(h)(x)&=\int_0^1h(t)D_n(x\dot+t)\,dt
 =w_n(x)\sum_{l=0}^{\infty}c_l(n)\int_0^1h(t)w_n(t)D_{2^l}(x\dot+t)\,dt\\
 &=w_n(x)\sum_{l=0}^{\infty}c_l(n)s_{2^l}(hw_n)(x)
 =w_n(x)\sum_{l=0}^{\infty}c_l(n)E_l(hw_n)(x).
\end{align*}
In the same way, the first expansion in \eqref{e2.8} gives
$s_n(h)=w_n\sum_{l=0}^{\infty}n_l(s_{2^{l+1}}-s_{2^l})(hw_n)$, and
$s_{2^{l+1}}-s_{2^l}=E_{l+1}-E_l=d_{l+1}$. Thus, we have
\begin{align}\label{e2.9}
 s_n(h)=w_n\sum_{l=0}^{\infty}c_l(n)E_l(hw_n),\qquad
 s_n(h)=w_n\sum_{l=0}^{\infty}n_ld_{l+1}(hw_n).
\end{align}

For any $k\in\mathbb N$, the remainder $n\bmod2^k$ has the digits $n_l$ at the
positions $l<k$ and zero digits at the positions $l\geq k$. Its digit at
position $k-1$ is followed by a zero at position $k$, so the definition
of $V$ gives
\[
 V(n\bmod2^k)=n_0+\sum_{l=1}^{k-1}|n_l-n_{l-1}|+n_{k-1},
\]
and \begin{align}
\sum_{l=0}^{k-1}|c_l(n)|
 =V(n\bmod2^k)-n_{k-1}\leq V(n\bmod2^k);\label{e2.10}
\end{align}
the second formula follows from the first, since $|c_0(n)|=n_0$ and
$|c_l(n)|=|n_l-n_{l-1}|$ for $l\in\mathbb N$. Moreover, for any $k\in\mathbb N$, the first formula,
applied with $k+1$ and with $k$, gives
\begin{align}
 V(n\bmod2^{k+1})-V(n\bmod2^k)
  &=|n_k-n_{k-1}|+n_k-n_{k-1}\nonumber\\
  &=2n_k\bigl(1-n_{k-1}\bigr)\geq0.\label{e2.11}
\end{align}
The last expression equals $2$ if $n_k=1$ and $n_{k-1}=0$, and $0$ otherwise. For
$k=0$, formula \eqref{e2.11} holds if we set $n_{-1}\coloneqq 0$, since
$n\bmod1=0$ and $V(n\bmod2)=2n_0$. Thus $V(n\bmod2^k)$ cannot decrease as
$k$ increases, and since it equals $V(n)$ for large $k$,
\begin{align}\label{e2.12}
 V(n\bmod2^k)\leq V(n)\qquad(k\in\mathbb Z_+).
\end{align}

To construct functions with large partial sums in Lemma~\ref{l5.2},
we use the following $L_1$ lower bound for the Walsh--Dirichlet kernel.
\begin{lemma}\label{l2.1}
For $d\in\mathbb Z_+$, let
$
 b_d\coloneqq \sum_{\substack{0\leq i<d\\i\text{ even}}}2^i
    =\sum_{i=0}^{\lceil d/2\rceil-1}4^i.
$
Then
\[
 \left\|D_{b_d}\right\|_{L_1[0,1)}\geq d/3.
\]
\end{lemma}
\begin{proof}
For $d=0$ both sides are zero, since $b_0=0$ and $D_0=0$. Let $d\in\mathbb N$.
Writing the powers of two in decreasing order, we have
$
 b_d=\sum_{p=1}^{\lceil d/2\rceil}2^{2(\lceil d/2\rceil-p)}.
$
Hence, Fine's formula for the Lebesgue constants
$\|D_n\|_{L_1[0,1)}$ \cite[Theorem~IX\,(i), p.~391]{Fine} gives
\begin{align*}
 \|D_{b_d}\|_{L_1[0,1)}
 &=\lceil d/2\rceil-\sum_{1\leq p<r\leq\lceil d/2\rceil}
   2^{2(\lceil d/2\rceil-r)-2(\lceil d/2\rceil-p)}
 =\lceil d/2\rceil-\sum_{p=1}^{\lceil d/2\rceil-1}
   \sum_{r=p+1}^{\lceil d/2\rceil}4^{p-r}\\
 &\geq\lceil d/2\rceil-\sum_{p=1}^{\lceil d/2\rceil-1}
   \sum_{r=p+1}^{\infty}4^{p-r}
 =\lceil d/2\rceil-\frac{\lceil d/2\rceil-1}{3}
 =\frac{2\lceil d/2\rceil+1}{3}\geq\frac d3.
\end{align*}
This finishes the proof of Lemma \ref{l2.1}.
\end{proof}

The analytic estimate needed in the sufficient direction is the theorem
of Billard \cite{Billard} and Sj\"olin \cite{Sjolin}. The \emph{maximal
partial sum operator} $s_*(h)\coloneqq \sup_{n\in\mathbb N}|s_n(h)|$ satisfies
\begin{align}\label{e2.13}
 \|s_*(h)\|_{L_2[0,1)}\leq B\|h\|_{L_2[0,1)}
 \qquad(h\in L_2[0,1)),
\end{align}
where $B\geq1$ is an absolute constant; see also
\cite{Thiele00,Demeter}.

\subsection{Weak type constants and the maximal principle}
For an operator $\Phi$ defined on $L_1[0,1)$, let
\begin{align}\label{e2.14}
 \|\Phi\|_{(1,1)}\coloneqq \sup_{\substack{h\in L_1[0,1),\ \|h\|_{L_1[0,1)}>0\\\varrho>0}}
       \frac{\varrho\,\lambda\{|\Phi(h)|>\varrho\}}{\|h\|_{L_1[0,1)}}
       \in[0,\infty],
\end{align}
the \emph{weak type $(1,1)$ constant} of $\Phi$. 
For any given set $F\subset\mathbb Z_+$, we write
\begin{align}\label{e2.15}
 s_*^F(f)\coloneqq \sup_{n\in F}|s_n(f)|,
\end{align}
and $s_*^\varnothing(f)\coloneqq 0$.

\begin{lemma}\label{l2.2}
Let $(\nu_j)_{j\in\mathbb N}$ be a sequence in $\mathbb N$ with
$\nu_1<\nu_2<\cdots$, and $F\coloneqq \{\nu_j:j\in\mathbb N\}$. The following
conditions are equivalent.
\begin{enumerate}[label=\textup{(\alph*)},leftmargin=*]
\item $\lim_{j\to\infty}s_{\nu_j}(f)=f$ almost everywhere for every $f\in L_1[0,1)$;
\item $s_*^F(f)<\infty$ almost everywhere for every $f\in L_1[0,1)$;
\item $\|s_*^F\|_{(1,1)}<\infty$.
\end{enumerate}
\end{lemma}
\begin{proof}
If (c) holds, then (a) follows by the usual density argument,
because the Walsh polynomials form a dense subset of $L_1[0,1)$ and
$s_{\nu_j}(P)=P$ for every Walsh polynomial $P$ and all sufficiently
large $j$. The implication (a)$\Rightarrow$(b) is obvious. The operators
$s_n$ are bounded on $L_1[0,1)$, and by \eqref{e2.2} they commute with
the dyadic translations. Hence (b)$\Rightarrow$(c) is the case $p=1$ of
the maximal principle of Stein \cite[Corollary~1, p.~148]{Stein} on the
dyadic group; see also \cite{Sawyer,On23}.
\end{proof}

\section{The main results}\label{s3}
In this section, we state a necessary and sufficient condition for the
selected Walsh--Fourier partial sums of $f$ to converge to $f$ almost
everywhere for every $f\in L_1[0,1)$. Upper and lower bounds for the weak
type $(1,1)$ constant of the associated maximal operator yield this
criterion. If the criterion fails, then there is an integrable
function whose selected partial sums are unbounded at every point.

\subsection{The criterion}
Let $F\subset\mathbb Z_+$. For $n\in F$ and $k\in\mathbb Z_+$ we write
\[
 n=2^kq+a,\qquad q\coloneqq \lfloor n/2^k\rfloor,\qquad a\coloneqq n\bmod2^k\in\{0,1,\ldots,2^k-1\}.
\]
The \emph{remainder} $a$ contains the first $k$ binary digits of $n$, and the
\emph{quotient} $q$ contains the remaining digits. The following definition
gives, at each scale, the set of distinct positive quotients associated
with remainders whose variation exceeds a fixed integer, and takes the
supremum of their number over all scales.
\begin{definition}\label{d3.1}
For any $k\in\mathbb Z_+$ and $K\geq1$, define
\begin{align*}
\mathcal H_{k,K}(F)
 &\coloneqq \bigl\{q\in\mathbb N:\ \exists\, a\in\{1,\ldots,2^k-1\},\quad
                   2^kq+a\in F,\quad V(a)>K\bigr\},\\
 \mathcal W_K(F)&\coloneqq \sup_{k\in\mathbb Z_+}\#\mathcal H_{k,K}(F).
\end{align*}
The value $\mathcal W_K(F)$ may be infinite.
\end{definition}

If $K\leq K'$, then
$\mathcal H_{k,K'}(F)\subset\mathcal H_{k,K}(F)$ for every $k$, and hence
$\mathcal W_{K'}(F)\leq\mathcal W_K(F)$. For example, if
$F=\{2^j:j\in\mathbb N\}$, then $\mathcal H_{k,K}(F)=\varnothing$ for all $k$ and
$K$. Hence $\mathcal W_K(F)=0$.

Our main result is the following theorem. The weak type constant
$\|s_*^F\|_{(1,1)}$ is defined in \eqref{e2.14} and \eqref{e2.15}.
\begin{theorem}\label{t3.2}
Let $(\nu_j)_{j\in\mathbb N}$ be a sequence in $\mathbb N$ with
$\nu_1<\nu_2<\cdots$, and put $F\coloneqq \{\nu_j:j\in\mathbb N\}$. The following
conditions are equivalent.
\begin{enumerate}[label=\textup{(\roman*)},leftmargin=*]
\item $\lim_{j\to\infty}s_{\nu_j}(f)=f$ almost everywhere for every $f\in L_1[0,1)$;
\item $\|s_*^F\|_{(1,1)}<\infty$;
\item $\mathcal W_K(F)<\infty$ for some $K\in\mathbb N$.
\end{enumerate}
\end{theorem}

\Cref{t3.2} answers Konyagin's problem \cite[Problem~3.3, p.~1398]{K06}
(see \cref{pr1.1}). For an arbitrary
sequence $(\nu_j)_{j\in\mathbb N}$ of positive integers, testing
convergence on each $w_m$ forces $\nu_j\to\infty$, since
$s_{\nu_j}(w_m)=0$ for $\nu_j\leq m$ and $|w_m|=1$. If $\nu_j\to\infty$,
pointwise convergence along $(\nu_j)_{j\in\mathbb N}$ is equivalent to
pointwise convergence along the increasing enumeration of
$F\coloneqq\{\nu_j:j\in\mathbb N\}$. Thus the required condition is
$\nu_j\to\infty$ and $\mathcal W_K(F)<\infty$ for some $K\in\mathbb N$.

The proof of \cref{t3.2} follows the statement of
\cref{t3.4} below and consists of three parts. The
equivalence (i)$\Leftrightarrow$(ii) is the translation principle of
\cref{l2.2}. The new part of the theorem is the equivalence of the
weak type estimate in (ii) and the counting condition in (iii). The
implication (iii)$\Rightarrow$(ii) follows from the upper estimate in
\cref{t3.4}\,(i) below, which is proved in \cref{s4}, and the
implication (ii)$\Rightarrow$(iii) follows from the lower estimate in
\cref{t3.4}\,(ii), which is proved in \cref{s5}.

The counting condition has an equivalent form that involves only one
integer. For $a\in\mathbb Z_+$ set $a_{-1}\coloneqq 0$ and define
\begin{align}\label{e3.1}
 \gamma(a)\coloneqq \sum_{l=0}^{\infty}a_l\bigl(1-a_{l-1}\bigr).
\end{align}
Since $|u-v|=2u(1-v)+v-u$ for $u,v\in\{0,1\}$, we have
\[
 V(a)=\sum_{l=0}^{\infty}|a_l-a_{l-1}|
 =2\gamma(a)+\sum_{l=0}^{\infty}(a_{l-1}-a_l)
 =2\gamma(a).
\]
The last sum is zero because $a_{-1}=0$ and $a_l=0$ for all sufficiently
large $l$.
\begin{corollary}\label{c3.3}
Let $(\nu_j)_{j\in\mathbb N}$ be a sequence in $\mathbb N$ with
$\nu_1<\nu_2<\cdots$, and put $F\coloneqq\{\nu_j:j\in\mathbb N\}$.
Then $s_{\nu_j}(f)\to f$ almost everywhere for every $f\in L_1[0,1)$
if and only if there exists $N\in\mathbb N$ such that
\[
 \#\left\{
 \left\lfloor\frac{n}{2^k}\right\rfloor:\
 n\in F,\ n\geq2^k,\ \gamma(n\bmod2^k)>N
 \right\}\leq N
 \qquad\text{for every } k\in\mathbb Z_+.
\]
Equivalently, $\mathcal W_{2N}(F)\leq N$ for some $N\in\mathbb N$.
\end{corollary}
\begin{proof}
For any $n=2^kq+a$ with $0\leq a<2^k$, we have $n\geq2^k$ if and
only if $q\geq1$. Since $V(a)=2\gamma(a)$, the inequality
$\gamma(a)>N$ is equivalent to $V(a)>2N$ and implies $a\ne0$.
Thus $\mathcal H_{k,2N}(F)$ consists exactly of the positive quotients
$\lfloor n/2^k\rfloor$ with $n\in F$ and $\gamma(n\bmod2^k)>N$.
The bounds $\#\mathcal H_{k,2N}(F)\leq N$ for all $k\in\mathbb Z_+$
are equivalent to $\mathcal W_{2N}(F)\leq N$.

If $\mathcal W_K(F)<\infty$ for some $K\in\mathbb N$, choose
$N\in\mathbb N$ with $N\geq\max\{\lceil K/2\rceil,\mathcal W_K(F)\}$.
Since $2N\geq K$, we have
\[
 \mathcal W_{2N}(F)\leq\mathcal W_K(F)\leq N.
\]
Conversely, if $\mathcal W_{2N}(F)\leq N$ for some $N\in\mathbb N$,
then $\mathcal W_K(F)<\infty$ with $K=2N$. The desired conclusion follows
from \cref{t3.2}.
\end{proof}

\subsection{The weak type constant}
In this subsection, we give upper and lower bounds for
$\|s_*^F\|_{(1,1)}$. The upper bound uses $\mathcal W_K(F)$, the
supremum of $\#\mathcal H_{k,K}(F)$ over all scales $k$. The lower
bound uses the size of a linearly independent subset of
$\mathcal H_{k,K}(F)$ at a single scale.

Let $\mathbb F_2\coloneqq\{0,1\}$ be the field with arithmetic modulo
two. We identify each nonnegative integer with the vector of its binary
digits over $\mathbb F_2$. A finite set $B\subseteq\mathbb Z_+$ is
\emph{linearly independent} if
\[
 \bigoplus_{q\in C}q\ne0
 \qquad\text{for every nonempty } C\subseteq B,
\]
where $\bigoplus$ denotes coordinatewise addition modulo two.
For any $A\subseteq\mathbb Z_+$, define its \emph{dimension} by
\[
 \dim_{\mathbb F_2}A\coloneqq
 \sup\{\#B:\ B\subseteq A,\ B\text{ finite and linearly independent}\}.
\]
For any $F\subseteq\mathbb Z_+$ and $K\geq1$, set
\begin{align}\label{e3.2}
 \operatorname{rank}_K(F)\coloneqq \sup_{k\in\mathbb Z_+}
       \dim_{\mathbb F_2}\mathcal H_{k,K}(F).
\end{align}
For $A\coloneqq\mathcal H_{k,K}(F)$, we have $\dim_{\mathbb F_2}A\leq\#A$.
If $r\coloneqq\dim_{\mathbb F_2}A<\infty$, choose a linearly independent
set $B\subseteq A$ with $\#B=r$. For any $q\in A\setminus B$, the set
$B\cup\{q\}$ is linearly dependent. A nonempty subset whose sum is zero
must contain $q$, since $B$ is independent. Hence
$q=\bigoplus_{p\in C}p$ for some $C\subseteq B$. Since $q>0$,
$C\ne\varnothing$. For $q\in B$, take $C=\{q\}$. Thus
\[
 A\subseteq\left\{\bigoplus_{q\in C}q:\ \varnothing\ne C\subseteq B\right\},
 \qquad \#A\leq 2^r-1.
\]
Taking suprema over $k$ gives
\begin{align}\label{e3.3}
 \operatorname{rank}_K(F)\leq\mathcal W_K(F),\qquad
 \mathcal W_K(F)\leq2^{\operatorname{rank}_K(F)}-1
       \quad\text{if }\operatorname{rank}_K(F)<\infty.
\end{align}

\begin{theorem}\label{t3.4}
Let $F\subset\mathbb Z_+$.
\begin{enumerate}[label=\textup{(\roman*)},leftmargin=*]
\item For all $K\in\mathbb N$ and $N\in\mathbb Z_+$ with
$\mathcal W_K(F)\leq N$, we have
\begin{align}\label{e3.4}
 \|s_*^F\|_{(1,1)}\leq C\bigl(K+\sqrt N\bigr),
\end{align}
where $C>0$ is an absolute constant.
\item If $d,L\in\mathbb N$ and
$\dim_{\mathbb F_2}\mathcal H_{k,d}(F)\geq L$
at some scale $k\in\mathbb Z_+$, then
\begin{align}\label{e3.5}
 \|s_*^F\|_{(1,1)}\geq\frac d{12}\big(1-e^{-L/2^d}\big).
\end{align}
\end{enumerate}
\end{theorem}

Part (i) of \cref{t3.4} is proved in \cref{s4}, and part (ii) in \cref{s5}. 

We now combine the estimates in \cref{t3.4} with the maximal principle
in \cref{l2.2} to prove \cref{t3.2}.

\begin{proof}[Proof of \cref{t3.2}]
Let $(\nu_j)_{j\in\mathbb N}$ and $F=\{\nu_j:j\in\mathbb N\}$ be as in
\cref{t3.2}. The equivalence of \textup{(i)} and \textup{(ii)} follows
from the equivalence of \textup{(a)} and \textup{(c)} in \cref{l2.2}.

To prove \textup{(iii)}$\Rightarrow$\textup{(ii)}, assume that $\mathcal W_K(F)<\infty$ for some $K\in\mathbb N$. Then
$\mathcal W_K(F)\in\mathbb Z_+$, and \cref{t3.4}\,(i) with
$N=\mathcal W_K(F)$ gives
\[
 \|s_*^F\|_{(1,1)}\lesssim K+\sqrt{\mathcal W_K(F)}<\infty,
\]
which is \textup{(ii)}.

To prove \textup{(ii)}$\Rightarrow$\textup{(iii)}, assume that $\|s_*^F\|_{(1,1)}<\infty$, and choose $d\in\mathbb N$ with
\[
 d>12\|s_*^F\|_{(1,1)}.
\]
Suppose that $\operatorname{rank}_d(F)=\infty$. Then,
by \eqref{e3.2}, for every $L\in\mathbb N$ there is a scale $k$ for which
$\dim_{\mathbb F_2}\mathcal H_{k,d}(F)\geq L$,
and \cref{t3.4}\,(ii) gives
\[
 \|s_*^F\|_{(1,1)}\geq\frac d{12}\big(1-e^{-L/2^d}\big)\qquad(L\in\mathbb N).
\]
Letting $L\to\infty$,
we obtain
\[
 \|s_*^F\|_{(1,1)}\geq\frac d{12}>\|s_*^F\|_{(1,1)},
\]
a contradiction. Hence $\operatorname{rank}_d(F)<\infty$, and
\eqref{e3.3} gives
$\mathcal W_d(F)\leq2^{\operatorname{rank}_d(F)}-1<\infty$. This is (iii)
with $K=d$. This finishes the proof of \cref{t3.2}
\end{proof}

\subsection{Divergence}
If the criterion fails, then \cref{t3.2} gives an integrable function
whose partial sums along $(\nu_j)$ do not converge almost everywhere. The
next two results show that the divergence is in fact much stronger. They
strengthen \cref{t3.2}, but they are not used in its proof, and they
are proved in \cref{s6}. The first result gives a nonnegative function
whose partial sums are unbounded almost everywhere.
\begin{corollary}\label{c3.5}
Let $(\nu_j)_{j\in\mathbb N}$ be a sequence in $\mathbb N$ with
$\nu_1<\nu_2<\cdots$, and $F\coloneqq \{\nu_j:j\in\mathbb N\}$. If
$\mathcal W_K(F)=\infty$ for every $K\in\mathbb N$, then there exists a
nonnegative function $f\in L_1[0,1)$ with $\|f\|_{L_1[0,1)}=1$ such that
\[
 \limsup_{j\to\infty}|s_{\nu_j}(f)(x)|=\infty
 \quad\text{for almost every }x\in[0,1).
\]
\end{corollary}

Without the nonnegativity of $f$, \cref{c3.5} follows from
\cref{t3.2}, which gives $\|s_*^F\|_{(1,1)}=\infty$, and from Stein's
theorem \cite[Corollary~2, p.~148]{Stein}. The second result gives
unboundedness at every point.
\begin{theorem}\label{t3.6}
Let $(\nu_j)_{j\in\mathbb N}$ be a sequence in $\mathbb N$ with
$\nu_1<\nu_2<\cdots$, and $F\coloneqq \{\nu_j:j\in\mathbb N\}$. If
$\mathcal W_K(F)=\infty$ for every $K\in\mathbb N$, then there exists a
function $f\in L_1[0,1)$ such that
\[
 \limsup_{j\to\infty}|s_{\nu_j}(f)(x)|=\infty
 \quad\text{for every }x\in[0,1).
\]
\end{theorem}

Together with \cref{t3.2}, \cref{t3.6} gives an alternative for
every increasing sequence $(\nu_j)_{j\in\mathbb N}$ in $\mathbb N$: either
$\lim_{j\to\infty}s_{\nu_j}(f)=f$ almost everywhere for every
$f\in L_1[0,1)$, or there is a function $f\in L_1[0,1)$ with
$\limsup_{j\to\infty}|s_{\nu_j}(f)(x)|=\infty$ for every $x\in[0,1)$.

\section{The sufficient estimate}\label{s4}
In this section we prove \cref{t3.4}\,(i).
To this end, fix $K\in\mathbb N$ and $N\in\mathbb Z_+$ with $\mathcal W_K(F)\leq N$, a function
$f\in L_1[0,1)$, a level $\varrho>0$, and a real parameter $\tau>0$ with
$2\lfloor K/2\rfloor\tau<1$. We first show that
\begin{align}\label{e4.1}
 \lambda\{s_*^F(f)>\varrho\}
 \leq\left\{\frac1\tau+
       \frac{2B^2(N+1)\tau}{(1-2\lfloor K/2\rfloor\tau)^2}\right\}
                         \frac{\|f\|_{L_1[0,1)}}{\varrho},
\end{align} where $B$ is the absolute constant
in the $L_2$ estimate \eqref{e2.13}. 
If $\|f\|_{L_1[0,1)}>\tau\varrho$, then
\[
 \lambda\{s_*^F(f)>\varrho\}\leq1<\frac{\|f\|_{L_1[0,1)}}{\tau\varrho},
\]
and \eqref{e4.1} holds. Assume henceforth that
$\|f\|_{L_1[0,1)}\leq\tau\varrho$.

\subsection{The decomposition}\label{s4.1}
Let $\mathcal B$ be the collection of the maximal dyadic intervals $I$ with
$\langle|f|\rangle_I>\tau\varrho$, and put $k_I\coloneqq -\log_2|I|$
$(I\in\mathcal B)$. Two intersecting dyadic intervals are nested, and no
interval of $\mathcal B$ contains another one. Hence $\mathcal B$ consists
of pairwise disjoint intervals, and their union is the set
\[
 \Omega\coloneqq \{M(f)>\tau\varrho\}.
\]
If $I^+$ is the
dyadic interval of length $2^{1-k_I}$ containing $I$, then $|I^+|=2|I|$,
and $\langle|f|\rangle_{I^+}\leq\tau\varrho$ by the maximality of $I$. It
follows that
\begin{align}\label{e4.2}
 \tau\varrho<\langle|f|\rangle_I
    \leq\frac{|I^+|}{|I|}\langle|f|\rangle_{I^+}\leq2\tau\varrho
 \qquad(I\in\mathcal B).
\end{align}
By \eqref{e2.1}, $|f|\leq M(f)\leq\tau\varrho$ almost everywhere on
$[0,1)\setminus\Omega$.

We follow the multifrequency
Calder\'on--Zygmund decomposition of Nazarov, Oberlin and Thiele
\cite[Theorem~1.1]{NOT} and its Walsh versions \cite[Section~3]{OT} and
\cite[proof of Lemma~3.2]{DL}. For any
$I\in\mathcal B$ with $|I|=2^{-k}$, define the finite set
\[
 Q_I\coloneqq \{0\}\cup\mathcal H_{k,K}(F),\qquad \#Q_I\leq N+1.
\]
Define the coefficients $\alpha_{I,q}$ and the function $p_I$ by
\begin{align*}
 \alpha_{I,q}\coloneqq \frac1{|I|}\int_I f(t)w_{2^kq}(t)\,dt,
 \qquad
 p_I\coloneqq \mathbf 1_I\sum_{q\in Q_I}\alpha_{I,q}w_{2^kq}.
\end{align*}
For
each $q\in Q_I$, \eqref{e2.7} gives
\[
 |\alpha_{I,q}|\leq\langle|f|\rangle_I,\qquad
 \int_Ip_Iw_{2^kq}\,d\lambda
       =|I|\,\alpha_{I,q}=\int_Ifw_{2^kq}\,d\lambda.
\]
Orthogonality \eqref{e2.7} and the Cauchy--Schwarz inequality give
\begin{align*}
 \|p_I\|_{L_2[0,1)}^2
   &=|I|\sum_{q\in Q_I}|\alpha_{I,q}|^2
     \leq(N+1)|I|\,\langle|f|\rangle_I^2,\\
 \|p_I\|_{L_1[0,1)}
   &\leq|I|^{1/2}\|p_I\|_{L_2[0,1)}
     \leq\sqrt{N+1}\,|I|\,\langle|f|\rangle_I
     =\sqrt{N+1}\int_I|f|\,d\lambda.
\end{align*}
Set
\[
 b_I\coloneqq f\mathbf 1_I-p_I\quad(I\in\mathcal B),\qquad
 g\coloneqq f\mathbf 1_{[0,1)\setminus\Omega}+\sum_{I\in\mathcal B}p_I.
\]
The sums converge in $L_1[0,1)$, because
\begin{align}
 \sum_{I\in\mathcal B}\|p_I\|_{L_1[0,1)}&\leq\sqrt{N+1}\,\|f\|_{L_1[0,1)},\nonumber\\*
 \sum_{I\in\mathcal B}\|b_I\|_{L_1[0,1)}
     &\leq\sum_{I\in\mathcal B}\left(\int_I|f|\,d\lambda+\|p_I\|_{L_1[0,1)}\right)
      \leq\bigl(1+\sqrt{N+1}\bigr)\|f\|_{L_1[0,1)}.\label{e4.3}
\end{align}
Moreover,
\begin{align}\label{e4.4}
 f=g+\sum_{I\in\mathcal B}b_I,\qquad
 \int_Ib_Iw_{2^{k_I}q}\,d\lambda=0\quad(q\in Q_I);
\end{align}
in particular $\int_Ib_I\,d\lambda=0$. The sets carrying the different
terms of $g$ are disjoint. Hence, by $|f|\leq\tau\varrho$ almost everywhere
on $[0,1)\setminus\Omega$ and by \eqref{e4.2},
\begin{align}
 \|g\|_{L_2[0,1)}^2
 &=\int_{[0,1)\setminus\Omega}|f|^2\,d\lambda+\sum_{I\in\mathcal B}\|p_I\|_{L_2[0,1)}^2\nonumber\\
 &\leq\tau\varrho\int_{[0,1)\setminus\Omega}|f|\,d\lambda
             +(N+1)\sum_{I\in\mathcal B}\langle|f|\rangle_I\int_I|f|\,d\lambda\nonumber\\
 &\leq\tau\varrho\int_{[0,1)\setminus\Omega}|f|\,d\lambda
                +2(N+1)\tau\varrho\int_\Omega|f|\,d\lambda\nonumber\\
 &\leq2(N+1)\tau\varrho\,\|f\|_{L_1[0,1)}.\label{e4.5}
\end{align}

\subsection{The error bound and the maximal estimate}

\begin{lemma}\label{l4.1}
For every $x\in[0,1)\setminus\Omega$ and every $n\in F$,
\begin{align}\label{e4.6}
 \left|s_n\!\left(\sum_{I\in\mathcal B}b_I\right)(x)\right|
       \leq2\lfloor K/2\rfloor\,M(f\mathbf 1_\Omega)(x).
\end{align}
\end{lemma}
\begin{proof}
Fix $I\in\mathcal B$, write $k\coloneqq k_I$ and $n=2^kq+a$ with $0\leq a<2^k$,
and let $x_I$ be the left endpoint of $I$. For any integrable function $h$
that vanishes outside $I$,
\begin{align}\label{e4.7}
 s_n(h)(x)=w_{2^kq}(x)D_a(x\dot+x_I)
                  \int_I h(t)w_{2^kq}(t)\,dt
 \qquad(x\notin I).
\end{align}
Indeed, let $x\notin I$. For almost every $t\in I$, the first $k$ binary
digits of $x\dot+t$ are $|x_l-t_l|=|x_l-(x_I)_l|$, $0\leq l<k$, and they
are not all zero, since $x$ and $t$ differ in one of the first $k$ digits;
hence $x\dot+t\notin I_k$, and \eqref{e2.5} gives $D_{2^kq}(x\dot+t)=0$.
Moreover, $a<2^k$, so that $D_a$ is constant on the dyadic intervals of
length $2^{-k}$, and $D_a(x\dot+t)=D_a(x\dot+x_I)$. By \eqref{e2.4},
\begin{align*}
 s_n(h)(x)&=\int_Ih(t)\bigl[D_{2^kq}(x\dot+t)
   +w_{2^kq}(x\dot+t)D_a(x\dot+t)\bigr]\,dt\\
 &=w_{2^kq}(x)D_a(x\dot+x_I)\int_Ih(t)w_{2^kq}(t)\,dt,
\end{align*}
which is \eqref{e4.7}. If $a=0$, then $D_a=0$ and the expression
vanishes. By \eqref{e4.4} and \eqref{e2.7}, we know that
\[
 \begin{cases}
\int_I b_Iw_{2^kq}\,d\lambda=0,&q\in Q_I,\\[4pt]
\int_I p_Iw_{2^kq}\,d\lambda=0,&q\notin Q_I.
 \end{cases}
\]
Consequently, \eqref{e4.7} with $h=b_I$ in the first case, and with
$h=p_I$ in the second case, gives that, for any $x\notin I$,
\[
 \begin{cases}
 s_n(b_I)(x)=0,&\lfloor n/2^k\rfloor\in Q_I,\\
 s_n(b_I)(x)=s_n(f\mathbf 1_I)(x)-s_n(p_I)(x)=s_n(f\mathbf 1_I)(x),&\lfloor n/2^k\rfloor\notin Q_I.
 \end{cases}
\]
For a fixed $n\in F$, put
\[
 \mathcal B_n\coloneqq \{I\in\mathcal B:\ \lfloor n/2^{k_I}\rfloor\notin Q_I\}.
\]
Since $0\in Q_I$, every $I\in\mathcal B_n$ satisfies $2^{k_I}\leq n$.
The intervals in $\mathcal B_n$ are disjoint, so $\#\mathcal B_n\leq n$.

The bound $|s_n(h)(x)|\leq n\|h\|_{L_1[0,1)}$ and \eqref{e4.3} allow
termwise application of $s_n$. The two cases above give
\begin{align}\label{e4.8}
 s_n\!\left(\sum_{I\in\mathcal B}b_I\right)(x)
      =\sum_{I\in\mathcal B_n}s_n(f\mathbf 1_I)(x)
 \qquad(x\notin\Omega).
\end{align}
If $\mathcal B_n=\varnothing$, \eqref{e4.6} follows. Otherwise put
$k_*\coloneqq\max\{k_I:I\in\mathcal B_n\}$.
For $x\notin\Omega$ and $l\geq k_I$, the interval $I_l(x)$ is disjoint
from $I$, so $E_l(f\mathbf 1_Iw_n)(x)=0$. Thus \eqref{e4.8}, the first
identity in \eqref{e2.9}, the inequality
$|E_l(f\mathbf 1_Iw_n)|\leq E_l(|f|\mathbf 1_I)$, the linearity of $E_l$, and
$\sum_{I\in\mathcal B_n}\mathbf 1_I\leq\mathbf 1_\Omega$ give
\begin{align*}
 \left|s_n\!\left(\sum_{I\in\mathcal B}b_I\right)(x)\right|
 &=\left|w_n(x)\sum_{l=0}^{k_*-1}c_l(n)
       \sum_{\substack{I\in\mathcal B_n\\k_I>l}}E_l(f\mathbf 1_Iw_n)(x)\right|\\
 &\leq\sum_{l=0}^{k_*-1}|c_l(n)|\,
       E_l\!\left(|f|\sum_{\substack{I\in\mathcal B_n\\k_I>l}}
                            \mathbf 1_I\right)(x)\\
 &\leq M(f\mathbf 1_\Omega)(x)\sum_{l=0}^{k_*-1}|c_l(n)|.
\end{align*}
At scale $k_*$ the positive quotient $\lfloor n/2^{k_*}\rfloor$ is
outside $\mathcal H_{k_*,K}(F)$, since some $I\in\mathcal B_n$ has
$k_I=k_*$. If we had $V(n\bmod2^{k_*})>K$, then $n\bmod2^{k_*}\ne0$,
since $V(0)=0$, and
$n=2^{k_*}\lfloor n/2^{k_*}\rfloor+n\bmod2^{k_*}\in F$ would give
$\lfloor n/2^{k_*}\rfloor\in\mathcal H_{k_*,K}(F)$ by \cref{d3.1}. Hence
$V(n\bmod2^{k_*})\leq K$. This variation is
even, so \eqref{e2.10} gives
\[
 \sum_{l=0}^{k_*-1}|c_l(n)|\leq V(n\bmod2^{k_*})\leq2\lfloor K/2\rfloor.
\]
This proves \eqref{e4.6}.
\end{proof}

\begin{proof}[Proof of \eqref{e4.1}]
Since $\Omega=\{M(f)>\tau\varrho\}$, we have
$M(f\mathbf 1_\Omega)(x)\leq M(f)(x)\leq\tau\varrho$ for
$x\in[0,1)\setminus\Omega$. By
\eqref{e4.4} and \eqref{e4.6},
\[
 |s_n(f)(x)|\leq|s_n(g)(x)|+2\lfloor K/2\rfloor\tau\varrho
 \qquad(n\in F)
\]
for almost every $x\notin\Omega$. Since $2\lfloor K/2\rfloor\tau<1$, this
implies that, up to a null set,
\[
 \{s_*^F(f)>\varrho\}\setminus\Omega
       \subset\bigl\{s_*^F(g)>(1-2\lfloor K/2\rfloor\tau)\varrho\bigr\}.
\]
By \eqref{e2.1}, $\lambda(\Omega)\leq\|f\|_{L_1[0,1)}/(\tau\varrho)$. We
now use this, \eqref{e2.13}, \eqref{e4.5}, Chebyshev's
inequality, and the pointwise bound $s_*^F(g)\leq s_*(g)$. They give
\begin{align*}
 \lambda\{s_*^F(f)>\varrho\}
 &\leq\lambda(\Omega)+\lambda\bigl\{s_*^F(g)>(1-2\lfloor K/2\rfloor\tau)\varrho\bigr\}\\
 &\leq\frac{\|f\|_{L_1[0,1)}}{\tau\varrho}
       +\frac{B^2\|g\|_{L_2[0,1)}^2}{(1-2\lfloor K/2\rfloor\tau)^2\varrho^2}\\
 &\leq\left\{\frac1\tau+
       \frac{2B^2(N+1)\tau}{(1-2\lfloor K/2\rfloor\tau)^2}\right\}
                         \frac{\|f\|_{L_1[0,1)}}{\varrho},
\end{align*}
which is \eqref{e4.1}.
\end{proof}

\begin{proof}[Proof of \cref{t3.4}\,(i)] To prove \eqref{e3.4}, choose
\[
 \tau=\frac1{2\lfloor K/2\rfloor+\sqrt2 B\sqrt{N+1}}.
\]
Then $1-2\lfloor K/2\rfloor\tau=\sqrt2B\sqrt{N+1}\,\tau$, and therefore
\begin{align*}
 \frac1\tau=2\lfloor K/2\rfloor+\sqrt2 B\sqrt{N+1}
 =\frac{2B^2(N+1)\tau}{(1-2\lfloor K/2\rfloor\tau)^2},
\end{align*}
so that the coefficient in \eqref{e4.1} equals
\[
 4\lfloor K/2\rfloor+2\sqrt2B\sqrt{N+1}
 \leq2K+2\sqrt2B\bigl(1+\sqrt N\bigr)\lesssim K+\sqrt N,
\]
since $K\geq1$. This proves \eqref{e3.4}.
\end{proof}

\section{The necessary estimate}\label{s5}
In this section we prove \cref{t3.4}\,(ii).

\subsection{Comparison with complete blocks}

For $Q\subset\mathbb Z_+$ and $k\in\mathbb Z_+$ we write
$2^kQ\coloneqq \{2^kq:q\in Q\}$ and $Q+1\coloneqq \{q+1:q\in Q\}$.
\begin{lemma}\label{l5.1}
For $Q\subset\mathbb Z_+$ and $k\in\mathbb Z_+$,
\begin{align}\label{e5.1}
 \left\|s_*^{2^kQ}\right\|_{(1,1)}=\left\|s_*^Q\right\|_{(1,1)}.
\end{align}
Let $F\subset\mathbb Z_+$ and let $T$ be a finite or countable index set.
For each $t\in T$, let $q_t\in\mathbb Z_+$ and $a_t\in\{1,\ldots,2^k-1\}$.
If $\{2^kq_t+a_t:t\in T\}\subset F$ and $Q=\{q_t:t\in T\}$, then
\begin{align}\label{e5.2}
 \left\|s_*^{Q+1}\right\|_{(1,1)}\leq\left\|s_*^F\right\|_{(1,1)},\qquad
 \left\|s_*^{2^kQ}\right\|_{(1,1)}\leq\left\|s_*^F\right\|_{(1,1)}.
\end{align}
\end{lemma}

\begin{proof}
Put $\pi_k(x)\coloneqq 2^kx-\lfloor 2^kx\rfloor$ for $x\in[0,1)$.
For $f\in L_1[0,1)$ and $0\leq i<2^k$, set
\[
 I_i\coloneqq[i2^{-k},(i+1)2^{-k}),\qquad
 f_i(v)\coloneqq f\bigl(2^{-k}(i+v)\bigr)\quad(v\in[0,1)).
\]
By \eqref{e2.3} and \eqref{e2.4}, the kernel $D_{2^kq}(x\dot+t)$ vanishes
for almost every $t\notin I_i$ when $x\in I_i$. On $I_i$, we use
\eqref{e2.6} and the substitution $v=2^kt-i$ to obtain
\begin{align}\label{e5.3}
 s_{2^kq}(f)(x)
 &=2^k\sum_{j=0}^{q-1}w_j(\pi_k(x))
     \int_{I_i} f(t)w_j(\pi_k(t))\,dt
 =s_q(f_i)(\pi_k(x)).
\end{align}
For $\varrho>0$, this gives
\begin{align*}
 \lambda\{s_*^{2^kQ}(f)>\varrho\}
 &=\sum_{i=0}^{2^k-1}2^{-k}\lambda\{s_*^Q(f_i)>\varrho\}\\
 &\leq\frac{\|s_*^Q\|_{(1,1)}}{\varrho}
        \sum_{i=0}^{2^k-1}2^{-k}\|f_i\|_{L_1[0,1)}
 =\frac{\|s_*^Q\|_{(1,1)}}{\varrho}\|f\|_{L_1[0,1)}
\end{align*}
when $\|s_*^Q\|_{(1,1)}<\infty$. Hence
$\|s_*^{2^kQ}\|_{(1,1)}\leq\|s_*^Q\|_{(1,1)}$, which also holds if the
right-hand side is infinite.

For $h\in L_1[0,1)$ and $0\leq r<2^k$, set
$h_{k,r}\coloneqq w_r(h\circ\pi_k)$. For $j\in\mathbb Z_+$ and $0\leq m<2^k$,
the product $w_rw_m$ is constant on each $I_i$. Using \eqref{e2.6}
and integrating first on each $I_i$, we get
\[
 \widehat{h_{k,r}}(2^kj+m)
 =\widehat h(j)\sum_{i=0}^{2^k-1}2^{-k}w_r(i2^{-k})w_m(i2^{-k})
 =\widehat h(j)\mathbf 1_{\{r=m\}} .
\]
Thus, for $q\in\mathbb Z_+$ and $0\leq a<2^k$,
\[
 s_{2^kq+a}(h_{k,r})
 =w_r\bigl(s_{q+\mathbf 1_{\{r<a\}}}(h)\circ\pi_k\bigr).
\]
The map $\pi_k$ preserves Lebesgue measure, and
$\|h_{k,r}\|_{L_1[0,1)}=\|h\|_{L_1[0,1)}$. Taking $a=r=0$ yields
$\|s_*^Q\|_{(1,1)}\leq\|s_*^{2^kQ}\|_{(1,1)}$, completing
\eqref{e5.1}. For the indices in \eqref{e5.2}, taking $r=0$ gives
\[
 s_*^{Q+1}(h)\circ\pi_k\leq s_*^F(h_{k,0}),
\]
while taking $r=2^k-1$ gives
\[
 s_*^Q(h)\circ\pi_k\leq s_*^F(h_{k,2^k-1}).
\]
The weak type estimates in \eqref{e5.2} follow from measure preservation
and \eqref{e5.1}.
\end{proof}

\subsection{Functions for the low remainders}
Goginava and Oniani showed that $|s_a(\operatorname{sgn}D_a)|$ is large
on a dyadic interval of length comparable to $1/a$
\cite[proof of Lemma~3.1]{GO}. In the following lemma, we use the sign
of $D_{b_d}$ and $d$ selected binary digits of $a$ to construct a bounded
function whose $a$th partial sum has absolute value
$\|D_{b_d}\|_{L_1[0,1)}$ on a set of measure $2^{-d}$.

\begin{lemma}\label{l5.2}
Let $k,d\in\mathbb N$ with $d\leq k$, and let $a\in\mathbb Z_+$ satisfy
$a<2^k$ and $V(a)>d$. Then there exist integers
$0\leq j_0<\cdots<j_{d-1}<k$ for which the map
\[
 \vartheta_a(x)\coloneqq \sum_{i=0}^{d-1}x_{j_i}2^{-i-1}
 \qquad(x\in[0,1)),
\]
satisfies $\lambda(\vartheta_a^{-1}(J))=2^{-d}$ for every dyadic interval
$J$ of length $2^{-d}$. Moreover, for any such $J$, there exists a real-valued function
$p_{a,J}$ on $[0,1)$, constant on each dyadic interval of length $2^{-k}$,
such that $|p_{a,J}|\leq1$ and
\begin{align}\label{e5.4}
 |s_a(p_{a,J})(x)|=\|D_{b_d}\|_{L_1[0,1)}\geq d/3
 \qquad\bigl(x\in\vartheta_a^{-1}(J)\bigr).
\end{align}
\end{lemma}
\begin{proof}
By \eqref{e3.1} and $V(a)=2\gamma(a)$, the first $k$ binary digits of
$a$ contain $V(a)/2$ runs of ones. Choose one digit from each run and
one from each gap between successive runs. In increasing order of their
indices, these $V(a)-1$ digits alternate, starting with $1$.
Since $V(a)-1\geq d$, let $j_0<\cdots<j_{d-1}$ be the indices of the
first $d$ chosen digits. By the definition of $b_d$ in \cref{l2.1},
\[
 a_{j_i}=(b_d)_i\qquad(0\leq i<d).
\]
For each dyadic interval $J$ of length $2^{-d}$, fixing the selected
digits leaves $k-d$ digits free among the first $k$. Thus
\begin{align}\label{e5.5}
 \lambda(\vartheta_a^{-1}(J))=2^{k-d}2^{-k}=2^{-d}.
\end{align}

For $0\leq m<2^d$ put $\iota(m)\coloneqq\sum_{i=0}^{d-1}m_i2^{j_i}$ and
set $a^\perp\coloneqq a-\iota(b_d)$. The selected digits of $a^\perp$
are all zero. By the definition of $\vartheta_a$ and the Walsh product rule,
\begin{align}\label{e5.6}
 w_{a^\perp}(w_m\circ\vartheta_a)=w_{a^\perp+\iota(m)}
 \qquad(0\leq m<2^d).
\end{align}
The map $\iota$ is strictly increasing because the highest bit at which
two integers differ is sent to the highest selected bit at which their
images differ. Hence $a^\perp+\iota(m)<a$ if and only if $m<b_d$.
Expanding any $h$ constant on dyadic intervals of length $2^{-d}$ in
the Walsh basis and retaining the frequencies less than $a$, we obtain
\begin{align}\label{e5.7}
 s_a\bigl(w_{a^\perp}(h\circ\vartheta_a)\bigr)
 =w_{a^\perp}\bigl(s_{b_d}(h)\circ\vartheta_a\bigr).
\end{align}

Let $y_J$ be the left endpoint of $J$, and set $\operatorname{sgn}(0)=0$.
Define
\begin{align}\label{e5.8}
 h_J(z)\coloneqq\operatorname{sgn}\bigl(D_{b_d}(z\dot+y_J)\bigr),
 \qquad p_{a,J}\coloneqq w_{a^\perp}(h_J\circ\vartheta_a).
\end{align}
The function $h_J$ is constant on the dyadic intervals of length
$2^{-d}$ since $b_d<2^d$ and $y_J$ is dyadic. Thus $p_{a,J}$ is real-valued,
constant on the dyadic intervals of length $2^{-k}$, and bounded in
absolute value by one. Since $s_{b_d}(h_J)$ is constant on $J$,
\eqref{e5.7} and the kernel formula give, for $x\in\vartheta_a^{-1}(J)$,
\begin{align*}
 |s_a(p_{a,J})(x)|
 &=|s_{b_d}(h_J)(y_J)|
 =\int_0^1|D_{b_d}(t\dot+y_J)|\,dt
 =\|D_{b_d}\|_{L_1[0,1)}\geq d/3.
\end{align*}
Here the last equality follows from the measure preservation of
$t\mapsto t\dot+y_J$, and the lower bound follows from \cref{l2.1}.
\end{proof}

\subsection{The product construction and the lower estimate}\label{s5.3}
For any finite family $\mathcal R_1,\ldots,\mathcal R_L$ of linear
operators defined on $L_1[0,1)$ and any $h\in L_1[0,1)$, let
\[
 \mathcal R^*(h)\coloneqq \max_{1\leq t\leq L}|\mathcal R_t(h)|.
\]

For related divergence constructions, see Schipp \cite{Schipp67,Schipp69}
and Konyagin \cite{K05}. Goginava and Oniani treat sequences with nested
spectra \cite[proof of Lemma~3.1]{GO}. In the following lemma, we use
the functions from \cref{l5.2} to form a product for a family of quotients
$q_1,\ldots,q_L$ that are linearly independent over $\mathbb F_2$. Each
factor uses the binary digits selected for its own remainder $a_t$.
This gives the lower bound \eqref{e5.9} in terms of $d$ and $L$.

\begin{lemma}\label{l5.3}
Let $k,d,L\in\mathbb N$ with $d\leq k$, and let $q_1,\ldots,q_L\in\mathbb N$
be linearly independent as binary vectors over
$\mathbb F_2$. Suppose that $a_t\in\mathbb Z_+$, $a_t<2^k$ and $V(a_t)>d$
for $1\leq t\leq L$. Then the operators
\[
 \mathcal R_t\coloneqq s_{2^kq_t+a_t}-s_{2^kq_t}\qquad(1\leq t\leq L)
\]
satisfy
\begin{align}\label{e5.9}
 \|\mathcal R^*\|_{(1,1)}\geq\|D_{b_d}\|_{L_1[0,1)}\big[1-(1-2^{-d})^L\big].
\end{align}
\end{lemma}
\begin{proof}
For each $t$, apply \cref{l5.2} with $a=a_t$; this is allowed, since
$d\leq k$ and $V(a_t)>d$. Let $\mathcal J$ be the set of all
$L$-tuples $(J_1,\ldots,J_L)$ in which each $J_t$ is one of the
$2^d$ dyadic intervals of length $2^{-d}$; the intervals $J_1,\ldots,J_L$
need not be different. Thus $\#\mathcal J=(2^d)^L=2^{dL}$. For
$(J_1,\ldots,J_L)\in\mathcal J$ put
\[
 U(J_1,\ldots,J_L)\coloneqq \bigcup_{t=1}^L\{x\in[0,1):\ \vartheta_{a_t}(x)\in J_t\}.
\]
Fix $x\in[0,1)$. For each $t$, exactly one of the $2^d$ intervals contains
$\vartheta_{a_t}(x)$, so that $J_t$ can be chosen in $2^d-1$ ways
with $\vartheta_{a_t}(x)\notin J_t$. Hence exactly $(2^d-1)^L$ tuples in
$\mathcal J$ satisfy $x\notin U(J_1,\ldots,J_L)$.
Interchanging the sum and the integral, we obtain
\begin{align*}
 \frac1{\#\mathcal J}\sum_{(J_1,\ldots,J_L)\in\mathcal J}
 \lambda\bigl(U(J_1,\ldots,J_L)\bigr)
 &=\int_0^1\frac{\#\{(J_1,\ldots,J_L)\in\mathcal J:\ x\in U(J_1,\ldots,J_L)\}}
   {\#\mathcal J}\,dx\\
 &=\frac{2^{dL}-(2^d-1)^L}{2^{dL}}=1-(1-2^{-d})^L.
\end{align*}
Fix a tuple $(J_1,\ldots,J_L)$ for which $U\coloneqq U(J_1,\ldots,J_L)$ has measure
at least this average, and set $p_t\coloneqq p_{a_t,J_t}$. Then, by \eqref{e5.4},
\begin{align}\label{e5.10}
 \lambda(U)\geq1-(1-2^{-d})^L,\quad |p_t|\leq1,\quad
 |s_{a_t}(p_t)(x)|=\|D_{b_d}\|_{L_1[0,1)}\quad(\vartheta_{a_t}(x)\in J_t).
\end{align}

Let
\[
 \mathcal P\coloneqq \prod_{t=1}^L\bigl(1+w_{2^kq_t}p_t\bigr).
\]
By \eqref{e5.10}, $1+w_{2^kq_t}p_t\geq1-|p_t|\geq0$ for every $t$.
Expanding the product and using $w_mw_n=w_{m\oplus n}$, we obtain
\begin{align*}
 \mathcal P
 &=\sum_{\Theta\subseteq\{1,\ldots,L\}}
   \prod_{\theta\in\Theta}\bigl(p_\theta w_{2^kq_\theta}\bigr)\\
 &=\sum_{\Theta\subseteq\{1,\ldots,L\}}
   \left(\prod_{\theta\in\Theta}p_\theta\right)
   w_{2^k\bigoplus_{\theta\in\Theta}q_\theta}.
\end{align*}
The term with $\Theta=\varnothing$ is $1$, with the convention
$\bigoplus_{\theta\in\varnothing}q_\theta\coloneqq0$.
Binary linear independence says exactly that, for $1\leq t\leq L$,
\[
 \bigoplus_{\theta\in\Theta}q_\theta=0\ \Longleftrightarrow\ \Theta=\varnothing,\qquad
 \bigoplus_{\theta\in\Theta}q_\theta=q_t\ \Longleftrightarrow\ \Theta=\{t\}.
\]
The functions $p_t$ are constant on the dyadic intervals of
length $2^{-k}$, while $E_k(w_{2^kq})=0$ for $q\in\mathbb N$ by \eqref{e2.7}.
Moreover, $w_{2^k\bigoplus_{\theta\in\Theta}q_\theta}w_{2^kq_t}
=w_{2^k((\bigoplus_{\theta\in\Theta}q_\theta)\oplus q_t)}$. Hence the first
equivalence leaves only the term $\Theta=\varnothing$ in $E_k(\mathcal P)$,
and the second leaves only the term $\Theta=\{t\}$ in
$E_k(\mathcal Pw_{2^kq_t})$:
\begin{align}\label{e5.11}
 \mathcal P\geq0,\qquad
 E_k(\mathcal P)=1,\qquad
 E_k\bigl(\mathcal Pw_{2^kq_t}\bigr)=p_t\qquad(1\leq t\leq L).
\end{align}
In particular,
$\|\mathcal P\|_{L_1[0,1)}=\int_0^1E_k(\mathcal P)\,d\lambda=1$.

For any $1\leq t\leq L$ and $0\leq m<2^k$, we have
$2^kq_t+m=(2^kq_t)\oplus m$. Since $w_m$ is constant on the dyadic
intervals of length $2^{-k}$, using \eqref{e5.11}, we obtain
\begin{align*}
 \widehat{\mathcal P}(2^kq_t+m)
 &=\int_0^1\mathcal Pw_{2^kq_t}w_m\,d\lambda
 =\int_0^1E_k\bigl(\mathcal Pw_{2^kq_t}\bigr)w_m\,d\lambda\\
 &=\int_0^1p_tw_m\,d\lambda=\widehat{p_t}(m).
\end{align*}
By the definition of $\mathcal R_t$ and $a_t<2^k$, it follows that
\begin{align*}
 \mathcal R_t(\mathcal P)
 &=\sum_{m=0}^{a_t-1}\widehat{\mathcal P}(2^kq_t+m)w_{2^kq_t}w_m\\
 &=w_{2^kq_t}\sum_{m=0}^{a_t-1}\widehat{p_t}(m)w_m
 =w_{2^kq_t}s_{a_t}(p_t).
\end{align*}
For any $x\in U$, choose $t$ such that $\vartheta_{a_t}(x)\in J_t$.
By \eqref{e5.10}, we obtain
\[
 \mathcal R^*(\mathcal P)(x)\geq|\mathcal R_t(\mathcal P)(x)|
 =|s_{a_t}(p_t)(x)|=\|D_{b_d}\|_{L_1[0,1)}.
\]
Therefore, for every level $0<\eta<\|D_{b_d}\|_{L_1[0,1)}$,
\[
 \lambda\{\mathcal R^*(\mathcal P)>\eta\}\geq\lambda(U)
       \geq1-(1-2^{-d})^L,\qquad \|\mathcal P\|_{L_1[0,1)}=1,
\]
and the definition \eqref{e2.14} gives
\[
 \|\mathcal R^*\|_{(1,1)}
 \geq\frac{\eta\,\lambda\{\mathcal R^*(\mathcal P)>\eta\}}{\|\mathcal P\|_{L_1[0,1)}}
 \geq\eta\big[1-(1-2^{-d})^L\big]\qquad(0<\eta<\|D_{b_d}\|_{L_1[0,1)}).
\]
Letting $\eta$ increase to $\|D_{b_d}\|_{L_1[0,1)}$ proves \eqref{e5.9}.
\end{proof}

\begin{proof}[Proof of \cref{t3.4}\,(ii)]
Suppose that $\dim_{\mathbb F_2}\mathcal H_{k,d}(F)\geq L$ for some $k$.
Choose linearly independent $q_1,\ldots,q_L\in\mathcal H_{k,d}(F)$ and
remainders $a_t$ with
\[
 1\leq a_t<2^k,\qquad V(a_t)>d,\qquad 2^kq_t+a_t\in F
 \qquad(1\leq t\leq L).
\]
Since $d<V(a_t)\leq k+1$, we have $d\leq k$. Let $Q=\{q_1,\ldots,q_L\}$.
For the operators $\mathcal R_t$ of \cref{l5.3},
\[
 \mathcal R^*(h)\leq s_*^F(h)+s_*^{2^kQ}(h).
\]
If $\|s_*^F\|_{(1,1)}=\infty$, the claimed bound is immediate.
Otherwise, by \cref{l5.1}, for every $\eta>0$,
\begin{align*}
 \lambda\{\mathcal R^*(h)>\eta\}
 &\leq\lambda\{s_*^F(h)>\eta/2\}+\lambda\{s_*^{2^kQ}(h)>\eta/2\}\\
 &\leq 4\|s_*^F\|_{(1,1)}\frac{\|h\|_{L_1[0,1)}}{\eta}.
\end{align*}
Combining $\|\mathcal R^*\|_{(1,1)}\leq4\|s_*^F\|_{(1,1)}$ with
\cref{l5.3}, we obtain
\[
 \|s_*^F\|_{(1,1)}
 \geq\tfrac14\|D_{b_d}\|_{L_1[0,1)}\big[1-(1-2^{-d})^L\big]
 \geq\frac d{12}\big(1-e^{-L/2^d}\big).
\]
The last inequality uses \cref{l2.1} and $(1-2^{-d})^L\leq e^{-L/2^d}$.
\end{proof}

\section{Divergence}\label{s6}
In this section, we prove \cref{c3.5,t3.6}. Throughout, $(\nu_j)_{j\in\mathbb N}$
and $F=\{\nu_j:j\in\mathbb N\}$ are as in these results.

\subsection{Divergence almost everywhere}
\begin{proof}[Proof of \cref{c3.5}]
By \cref{t3.2}, $\|s_*^F\|_{(1,1)}=\infty$, and \cref{l2.2} gives an
integrable function whose maximal function $s_*^F$ is infinite on a set of
positive measure. The maximal function of a sum is at most the sum of the
maximal functions of the summands. Taking the positive and the negative
part, we therefore obtain a nonnegative function $g\in L_1[0,1)$ such that
\[
 A\coloneqq \{x\in[0,1):\ s_*^F(g)(x)=\infty\},\qquad \lambda(A)>0;
\]
indeed, if both parts had maximal functions that are finite almost
everywhere, then so would the original function.

Since $\sum_{\ell=1}^{\infty}\lambda(A)=\infty$, the lemma of Stein
\cite[Lemma~1, p.~146]{Stein}, applied on the dyadic group to the sets
$A,A,\ldots$, gives points $t_1,t_2,\ldots\in[0,1)$ such that
\begin{align}\label{e6.1}
 \lambda\Bigl(\bigcup_{\ell=1}^{\infty}(A\dot+t_\ell)\Bigr)=1,\qquad
 A\dot+t\coloneqq \{x\dot+t:x\in A\}.
\end{align}

Let $\mathcal X\coloneqq \prod_{\ell=1}^{\infty}[0,2^{-\ell}]$, and let $\mu$ be the product of
the normalized Lebesgue measures $2^\ell\,d\xi_\ell$ on the intervals
$[0,2^{-\ell}]$. For $\xi=(\xi_\ell)_{\ell\in\mathbb N}\in\mathcal X$ set
\[
 f_\xi\coloneqq \sum_{\ell=1}^{\infty}\xi_\ell\, g(\,\cdot\,\dot+t_\ell).
\]
The series converges in $L_1[0,1)$ for every $\xi\in\mathcal X$ and defines a
nonnegative function, since
\[
 \sum_{\ell=1}^{\infty}
 \|\xi_\ell\, g(\,\cdot\,\dot+t_\ell)\|_{L_1[0,1)}
 =\|g\|_{L_1[0,1)}\sum_{\ell=1}^{\infty}\xi_\ell\leq\|g\|_{L_1[0,1)}.
\]
For every fixed $n$, the bound
$\|s_n(h)\|_{L_\infty[0,1)}\leq n\|h\|_{L_1[0,1)}$ permits the termwise
application of $s_n$, with uniform convergence of the resulting series.
By \eqref{e2.2}, we have, for almost every $x$,
\begin{align}\label{e6.2}
 s_{\nu_j}(f_\xi)(x)=\xi_\ell\,s_{\nu_j}(g)(x\dot+t_\ell)
 +\sum_{\substack{i=1\\i\ne\ell}}^{\infty}\xi_i\,s_{\nu_j}(g)(x\dot+t_i)
 \qquad(j,\ell\in\mathbb N).
\end{align}

Fix $\ell$ and a point $x\in A\dot+t_\ell$ for which \eqref{e6.2}
holds and $x\dot+t_\ell\in A$, and fix all coordinates of $\xi$ except
$\xi_\ell$. Then
$\sup_{j\in\mathbb N}|s_{\nu_j}(g)(x\dot+t_\ell)|=s_*^F(g)(x\dot+t_\ell)=\infty$.
There is at most one real number $u$ for which the sequence
\[
 \Bigl(u\,s_{\nu_j}(g)(x\dot+t_\ell)
 +\sum_{\substack{i=1\\i\ne\ell}}^{\infty}\xi_i\,s_{\nu_j}(g)(x\dot+t_i)\Bigr)_{j\in\mathbb N}
\]
is bounded. Indeed, if two distinct numbers $u$ and $u'$ had this
property, then the difference of the two sequences,
$\bigl((u-u')s_{\nu_j}(g)(x\dot+t_\ell)\bigr)_{j\in\mathbb N}$, would be bounded,
which contradicts $\sup_{j\in\mathbb N}|s_{\nu_j}(g)(x\dot+t_\ell)|=\infty$.
By \eqref{e6.2}, for $u=\xi_\ell$ this sequence is
$\bigl(s_{\nu_j}(f_\xi)(x)\bigr)_{j\in\mathbb N}$. Hence $s_*^F(f_\xi)(x)<\infty$
holds for at most one $\xi_\ell\in[0,2^{-\ell}]$, and a
single point has measure zero with respect to $2^\ell\,d\xi_\ell$. By
Fubini's theorem on $\mathcal X$, it follows that
\begin{align}\label{e6.3}
 \mu\bigl\{\xi\in\mathcal X:\ s_*^F(f_\xi)(x)<\infty\bigr\}=0
 \qquad\text{for almost every }x\in A\dot+t_\ell.
\end{align}
For each fixed $n$, the function $(x,\xi)\mapsto s_n(f_\xi)(x)$ is the sum
of a uniformly convergent series of measurable functions of $(x,\xi)$, and
hence it is jointly measurable; so is $(x,\xi)\mapsto s_*^F(f_\xi)(x)$,
since $F$ is countable. By \eqref{e6.1}, the equality in
\eqref{e6.3} holds for almost every $x\in[0,1)$. Hence Fubini's theorem on $[0,1)\times\mathcal X$ gives
\[
 \int_{\mathcal X}\lambda\{x\in[0,1):\ s_*^F(f_\xi)(x)<\infty\}\,d\mu(\xi)
 =\int_0^1\mu\bigl\{\xi\in\mathcal X:\ s_*^F(f_\xi)(x)<\infty\bigr\}\,dx
 =0.
\]
Choose $\xi\in\mathcal X$ for which $s_*^F(f_\xi)=\infty$ almost everywhere. Then
$\|f_\xi\|_{L_1[0,1)}>0$. Each fixed partial sum is finite, so that the
unboundedness implies
\[
 \limsup_{j\to\infty}|s_{\nu_j}(f_\xi)(x)|=\infty
 \quad\text{for almost every }x\in[0,1).
\]
The function $f\coloneqq f_\xi/\|f_\xi\|_{L_1[0,1)}$ has all the required properties.
\end{proof}

\subsection{Divergence everywhere}
In this subsection, we prove \cref{t3.6} by adapting the construction
of Goginava and Oniani for nested spectra \cite[Section~4]{GO} and
Oniani's argument for passing from almost everywhere divergence to
divergence everywhere \cite[Lemmas~1 and~2]{On22}. At each stage,
the set on which all selected partial sums of a Walsh polynomial are
small is a union of dyadic intervals of the same length, so its
indicator is a Walsh polynomial. We multiply this indicator by a Walsh
function and a scalar to obtain a large difference of partial sums on
that set.

For any Walsh polynomial $P$ and any $Z\subseteq\mathbb Z_+$, we say
that the \emph{frequencies of $P$ lie in $Z$} if $\widehat P(m)=0$ for
every $m\in\mathbb Z_+\setminus Z$. By the definition of $s_n$, for
any $n,n'\in\mathbb Z_+$,
\begin{align*}
 s_n(P)&=P &&\text{if $n\geq n'$ and the frequencies of $P$ are below $n'$},\\
 s_n(P)&=0 &&\text{if $n\leq n'$ and the frequencies of $P$ are at least $n'$}.
\end{align*}

The following lemma uses the function from \cref{c3.5} to construct,
for any $j_0\in\mathbb Z_+$, a Walsh polynomial $P$ whose frequencies
are at least $2^{j_0}$. The maximum of $|s_n(P)|$ over finitely many
indices $n\in F$ is large outside a set of small measure.

\begin{lemma}\label{l6.1}
Let $(\nu_j)_{j\in\mathbb N}$ be a sequence in $\mathbb N$ with
$\nu_1<\nu_2<\cdots$, put $F\coloneqq \{\nu_j:j\in\mathbb N\}$, and assume that
$\mathcal W_K(F)=\infty$ for every $K\in\mathbb N$. For all
$\ell\in\mathbb N$ and $j_0\in\mathbb Z_+$ there exist a finite nonempty set
$F'\subset\{n\in F:\ n>2^{j_0}\}$ and a Walsh polynomial $P$ such
that
\begin{align}\label{e6.4}
\begin{split}
 &\|P\|_{L_1[0,1)}\leq4,\qquad
 \widehat P(m)=0\quad(0\leq m<2^{j_0}),\\
 &\lambda\Bigl\{x\in[0,1):\ \max_{n\in F'}|s_n(P)(x)|\leq4^\ell\Bigr\}\leq2^{1-\ell}.
\end{split}
\end{align}
\end{lemma}
\begin{proof}
Let $f_0$ be the nonnegative function of \cref{c3.5}. For almost every
$x$, the sequence $(s_{\nu_j}(f_0)(x))_{j\in\mathbb N}$ is unbounded, and
therefore
\[
 \sup_{n\in F,\ n>2^{j_0}}|s_n(f_0)(x)|=\infty
 \quad\text{for almost every }x\in[0,1).
\]
Consider $N\in\mathbb N$ with $F\cap(2^{j_0},N]\ne\varnothing$.

The sets
$\{x\in[0,1):\ \max_{n\in F,\ 2^{j_0}<n\leq N}|s_n(f_0)(x)|\leq2\cdot4^\ell+1\}$
decrease as $N$ increases, and their intersection is a null set. By the
continuity of measure, there is $N\in\mathbb N$ with $N>2^{j_0}$ and
\begin{align}\label{e6.5}
 \lambda\Bigl\{x\in[0,1):\ \max_{n\in F,\ 2^{j_0}<n\leq N}|s_n(f_0)(x)|\leq2\cdot4^\ell+1\Bigr\}
 \leq2^{-\ell}.
\end{align}
Let $F'\coloneqq\{n\in F:\ 2^{j_0}<n\leq N\}$. By density, choose a
Walsh polynomial $h$ with $\|f_0-h\|_{L_1[0,1)}\leq1/N$. Since
$\|f_0\|_{L_1[0,1)}=1$, we have $\|h\|_{L_1[0,1)}\leq2$ and
\[
 |s_n(f_0-h)|\leq n\|f_0-h\|_{L_1[0,1)}\leq1\qquad(n\in F').
\]

Let $P\coloneqq h-E_{j_0}(h)$. Since $E_{j_0}(h)=s_{2^{j_0}}(h)$ and
$E_{j_0}$ is contractive on $L_1[0,1)$, we obtain
\[
 \widehat P(m)=0\quad(0\leq m<2^{j_0}),\qquad
 \|P\|_{L_1[0,1)}\leq2\|h\|_{L_1[0,1)}\leq4.
\]
For any $n\in F'$, we have $n>2^{j_0}$ and hence
$s_n(P)=s_n(h)-E_{j_0}(h)$. Using $|E_{j_0}(h)|\leq M(h)$ and the
bound on $s_n(f_0-h)$, we obtain
\[
 \max_{n\in F'}|s_n(f_0)|\leq\max_{n\in F'}|s_n(P)|+1+M(h).
\]
Thus, by \eqref{e6.5} and \eqref{e2.1},
\begin{align*}
 \lambda\Bigl\{\max_{n\in F'}|s_n(P)|\leq4^\ell\Bigr\}
 &\leq\lambda\Bigl\{\max_{n\in F'}|s_n(f_0)|\leq2\cdot4^\ell+1\Bigr\}
   +\lambda\{M(h)>4^\ell\}\\
 &\leq2^{-\ell}+\frac{\|h\|_{L_1[0,1)}}{4^\ell}
 \leq2^{-\ell}+2\cdot4^{-\ell}\leq2^{1-\ell}.
\end{align*}
This proves \eqref{e6.4}.
\end{proof}

\begin{proof}[Proof of \cref{t3.6}]
Put $m''_0\coloneqq0$. At stage $\ell\in\mathbb N$, choose $m_\ell\in F$
with $m_\ell>m''_{\ell-1}$ and $j_\ell\in\mathbb Z_+$ with
$m_\ell\leq2^{j_\ell}$. Apply \cref{l6.1} with $\ell$ and $j_0=j_\ell$
to obtain
$F'_\ell$ and $P_\ell$, and set
\[
 Y_\ell\coloneqq
 \Bigl\{x\in[0,1):\ \max_{n\in F'_\ell}|s_n(P_\ell)(x)|\leq4^\ell\Bigr\},
 \qquad \lambda(Y_\ell)\leq2^{1-\ell}.
\]
Choose $k_\ell\in\mathbb N$ with $k_\ell>j_\ell$ so that $\max F'_\ell$
and all frequencies of $P_\ell$ are less than $2^{k_\ell}$. Since $F$ is unbounded, we can
choose $m'_\ell,m''_\ell\in F$ so that the construction satisfies
\begin{align}\label{e6.6}
\begin{gathered}
 m''_{\ell-1}<m_\ell\leq2^{j_\ell}<\min F'_\ell,\qquad
 \max F'_\ell<2^{k_\ell}\leq m'_\ell<m''_\ell,\\
 \widehat{P_\ell}(m)=0\quad(m\notin[2^{j_\ell},2^{k_\ell})),\qquad
 m''_\ell-m'_\ell\geq2^{k_\ell+1}.
\end{gathered}
\end{align}

For any $n\in F'_\ell$, we have $n<2^{k_\ell}$, so $s_n(P_\ell)$ is
constant on each dyadic interval of length $2^{-k_\ell}$. Thus $Y_\ell$
is a union of such intervals and
\[
 \mathbf 1_{Y_\ell}=E_{k_\ell}(\mathbf 1_{Y_\ell})
 =s_{2^{k_\ell}}(\mathbf 1_{Y_\ell}).
\]
Let
\[
 u_\ell\coloneqq2^{k_\ell}\left\lceil m'_\ell/2^{k_\ell}\right\rceil,
 \qquad
 G_\ell\coloneqq2^{\ell/2}w_{u_\ell}\mathbf 1_{Y_\ell}.
\]
Since $u_\ell$ is a multiple of $2^{k_\ell}$, the Walsh product rule gives
\[
 w_{u_\ell}w_m=w_{u_\ell\oplus m}=w_{u_\ell+m}
 \qquad(0\leq m<2^{k_\ell}).
\]
By the definition of $u_\ell$ and \eqref{e6.6}, we have
\[
 m'_\ell\leq u_\ell,\qquad
 u_\ell+2^{k_\ell}<m'_\ell+2^{k_\ell+1}\leq m''_\ell.
\]
Hence the frequencies of $G_\ell$ lie in
$[u_\ell,u_\ell+2^{k_\ell})\subseteq[m'_\ell,m''_\ell)$. Using
$|w_{u_\ell}|=1$ and $\lambda(Y_\ell)\leq2^{1-\ell}$, we obtain
\[
 |G_\ell|=2^{\ell/2}\mathbf 1_{Y_\ell},\qquad
 \|G_\ell\|_{L_1[0,1)}=2^{\ell/2}\lambda(Y_\ell)\leq2^{1-\ell/2}.
\]
The series
\[
 f\coloneqq\sum_{\ell=1}^\infty(2^{-\ell}P_\ell+G_\ell)
\]
converges in $L_1[0,1)$, since $\|P_\ell\|_{L_1[0,1)}\leq4$ and
$\sum_{\ell=1}^{\infty}2^{1-\ell/2}<\infty$.

For each $n$, we can apply $s_n$ termwise by
$|s_n(h)(x)|\leq n\|h\|_{L_1[0,1)}$. By \eqref{e6.6}, the frequencies
from stages before $\ell$ are below $m_\ell$, and those from stages
after $\ell$ are at least $m''_\ell$. Within stage $\ell$, the
frequencies of $P_\ell$ are below $m'_\ell$, while those of $G_\ell$
are in $[m'_\ell,m''_\ell)$. It follows that
\begin{align}\label{e6.7}
 s_n(f)-s_{m_\ell}(f)&=2^{-\ell}s_n(P_\ell)\qquad(n\in F'_\ell),\\
 s_{m''_\ell}(f)-s_{m'_\ell}(f)&=G_\ell.\label{e6.8}
\end{align}

Using \eqref{e6.7}, \eqref{e6.8}, and the definition of $Y_\ell$, we
obtain, for any $x\in[0,1)$,
\begin{align*}
 \max_{n\in F'_\ell}\bigl|s_n(f)(x)-s_{m_\ell}(f)(x)\bigr|
 &>2^\ell &&\text{if }x\notin Y_\ell,\\
 \bigl|s_{m''_\ell}(f)(x)-s_{m'_\ell}(f)(x)\bigr|
 &=|G_\ell(x)|=2^{\ell/2} &&\text{if }x\in Y_\ell.
\end{align*}
By the triangle inequality,
\[
 \max_{n\in F'_\ell\cup\{m_\ell,m'_\ell,m''_\ell\}}
             |s_n(f)(x)|\geq2^{\ell/2-1}\qquad(x\in[0,1)).
\]
By \eqref{e6.6},
\[
 F'_\ell\cup\{m_\ell,m'_\ell,m''_\ell\}\subseteq F\cap[m_\ell,\infty),
 \qquad m_\ell\to\infty.
\]
Letting $\ell\to\infty$, we obtain
$\limsup_{j\to\infty}|s_{\nu_j}(f)(x)|=\infty$ for every $x\in[0,1)$.
This completes the proof of \cref{t3.6}.
\end{proof}

\section{Consequences and examples}\label{s7}
In this section, we apply \cref{t3.2,t3.6} to recover earlier results
and to give two examples. In \cref{ex7.6}, we obtain almost everywhere
convergence to $f$ for every $f\in L_1[0,1)$ along an index set that
cannot be covered by finitely many sets with bounded variation or
separated spectra. In \cref{ex7.7}, we obtain everywhere unbounded
partial sums for some $f\in L_1[0,1)$ along a sequence with no
subsequence at bounded distance from a sequence with nested spectra.

\subsection{Recovery of earlier results}\label{s7.1}
In this subsection, we apply \cref{t3.2} to sequences with bounded
variation \cite{Fine,SWS} and to sequences with separated spectra
\cite{K93}. We then apply \cref{c3.5,t3.6} to sequences with nested
spectra and unbounded variation, as studied by Goginava and Oniani
\cite{GO}. We first recall the three conditions on a sequence.

\begin{definition}\label{d7.1}
Let $(\nu_j)_{j\in\mathbb N}$ be a sequence in $\mathbb N$ with
$\nu_1<\nu_2<\cdots$.
It has \emph{bounded variation} if $\sup_{j\in\mathbb N}V(\nu_j)<\infty$, it has
\emph{separated spectra} if
\begin{align}\label{e7.1}
 \max\operatorname{Sp}(\nu_j)<\min\operatorname{Sp}(\nu_{j+1})
 \qquad(j\in\mathbb N),
\end{align}
and it has \emph{nested spectra} if
\begin{align}\label{e7.2}
 \operatorname{Sp}(\nu_{j+1})\cap[0,\max\operatorname{Sp}(\nu_j)]
 =\operatorname{Sp}(\nu_j)\qquad(j\in\mathbb N).
\end{align}
A subset of $\mathbb N$ has one of these properties if its increasing
enumeration has it.
\end{definition}

The sequences with separated spectra and those with nested spectra
form disjoint classes, and each of these classes meets the class of
sequences with bounded variation, but neither contains it or is contained
in it. Indeed, by \eqref{e7.1} the left-hand side of \eqref{e7.2} is
empty, whereas $\operatorname{Sp}(\nu_j)\neq\varnothing$ because
$\nu_j\in\mathbb N$. The sequences $(2^j)_{j\in\mathbb N}$ and $(2^j-1)_{j\in\mathbb N}$ have
bounded variation; the first one has separated but not nested spectra, and
the second one has nested but not separated spectra. On the other hand, a
sequence with separated spectra can have unbounded variation: the sequence
$\nu_j\coloneqq \sum_{i=0}^{j-1}2^{j^2+2i}$ $(j\in\mathbb N)$ satisfies $V(\nu_j)=2j$ and
$\max\operatorname{Sp}(\nu_j)=j^2+2j-2<(j+1)^2=\min\operatorname{Sp}(\nu_{j+1})$.
Similarly, the sequence $\nu_j\coloneqq \sum_{i=0}^{j}4^i$ $(j\in\mathbb N)$ has nested
spectra and $V(\nu_j)=2j+2$.

\begin{corollary}[{\cite{Fine,SWS}}]\label{c7.2}
Let $(\nu_j)_{j\in\mathbb N}$ be an increasing sequence in $\mathbb N$
with bounded variation. Then $\lim_{j\to\infty}s_{\nu_j}(f)=f$ almost
everywhere for every $f\in L_1[0,1)$.
\end{corollary}
\begin{proof}
Put $F\coloneqq \{\nu_j:j\in\mathbb N\}$, and choose $K\in\mathbb N$ with
$V(\nu_j)\leq K$ for all $j\in\mathbb N$. By \eqref{e2.12},
\[
 V(\nu_j\bmod2^k)\leq V(\nu_j)\leq K\qquad(j\in\mathbb N,\ k\in\mathbb Z_+).
\]
Hence, by \cref{d3.1}, $\mathcal H_{k,K}(F)=\varnothing$ for every
$k\in\mathbb Z_+$, so that $\mathcal W_K(F)=0$, and \cref{t3.2} gives the
convergence.
\end{proof}

\begin{corollary}[{\cite{K93}}]\label{c7.3}
Let $(\nu_j)_{j\in\mathbb N}$ be an increasing sequence in $\mathbb N$
with separated spectra. Then $\lim_{j\to\infty}s_{\nu_j}(f)=f$ almost
everywhere for every $f\in L_1[0,1)$.
\end{corollary}
\begin{proof}
For $n\in\mathbb N$ and $k\in\mathbb Z_+$, the quotient $\lfloor n/2^k\rfloor$
contains the binary digits of $n$ at the positions at least $k$, and the
remainder $n\bmod2^k$ contains those at the positions below $k$. Hence
\begin{align}\label{e7.3}
 \lfloor n/2^k\rfloor\in\mathbb N\iff\max\operatorname{Sp}(n)\geq k,\qquad
 n\bmod2^k\in\mathbb N\iff\min\operatorname{Sp}(n)<k.
\end{align}
Put $F\coloneqq \{\nu_j:j\in\mathbb N\}$ and fix $k\in\mathbb Z_+$. By \cref{d3.1}, every
element of $\mathcal H_{k,1}(F)$ is the quotient
$\lfloor\nu_j/2^k\rfloor\in\mathbb N$ of a term $\nu_j$ with
$\nu_j\bmod2^k\in\mathbb N$, and for this $j$, \eqref{e7.3} gives
\[
 \min\operatorname{Sp}(\nu_j)<k\leq\max\operatorname{Sp}(\nu_j).
\]
By \eqref{e7.1}, the intervals
$[\min\operatorname{Sp}(\nu_j),\max\operatorname{Sp}(\nu_j)]$, $j\in\mathbb N$, are
pairwise disjoint. Hence at most one $j$ satisfies these inequalities,
and $\#\mathcal H_{k,1}(F)\leq1$. Thus $\mathcal W_1(F)\leq1$, and
\cref{t3.2} gives the convergence.
\end{proof}

\begin{lemma}\label{l7.4}
Let $(\nu_j)_{j\in\mathbb N}$ be an increasing sequence in $\mathbb N$
with nested spectra. Then
\begin{align}\label{e7.4}
 \nu_{j'}\bmod2^{1+\max\operatorname{Sp}(\nu_j)}=\nu_j\qquad(1\leq j\leq j').
\end{align}
In particular, every subsequence of $(\nu_j)_{j\in\mathbb N}$ has nested spectra.
\end{lemma}
\begin{proof}
By \eqref{e7.2}, $\operatorname{Sp}(\nu_i)\subset\operatorname{Sp}(\nu_{i+1})$
for $i\in\mathbb N$, and therefore
\[
 \max\operatorname{Sp}(\nu_j)\leq\max\operatorname{Sp}(\nu_i)
 \qquad(1\leq j\leq i).
\]
Let $1\leq j\leq i$. Then
$[0,\max\operatorname{Sp}(\nu_j)]\subset[0,\max\operatorname{Sp}(\nu_i)]$,
and the intersection of both sides of \eqref{e7.2}, taken with $i$ in
place of $j$, with $[0,\max\operatorname{Sp}(\nu_j)]$ gives
\begin{align}\label{e7.5}
 \operatorname{Sp}(\nu_{i+1})\cap[0,\max\operatorname{Sp}(\nu_j)]
 =\operatorname{Sp}(\nu_i)\cap[0,\max\operatorname{Sp}(\nu_j)]
 \qquad(1\leq j\leq i).
\end{align}
Applying \eqref{e7.5} with $i=j,j+1,\ldots,j'-1$, we obtain
\[
 \operatorname{Sp}(\nu_{j'})\cap[0,\max\operatorname{Sp}(\nu_j)]
 =\operatorname{Sp}(\nu_j)\cap[0,\max\operatorname{Sp}(\nu_j)]
 =\operatorname{Sp}(\nu_j)\qquad(1\leq j\leq j').
\]
Since $\operatorname{Sp}(n\bmod2^k)=\operatorname{Sp}(n)\cap[0,k-1]$ for
$n,k\in\mathbb Z_+$, this is \eqref{e7.4}. If $j_1<j_2<\cdots$, then
\eqref{e7.4} with $j=j_l$ and $j'=j_{l+1}$ gives \eqref{e7.2} for the
subsequence $(\nu_{j_l})_{l\in\mathbb N}$.
\end{proof}

\begin{corollary}\label{c7.5}
Let $(\nu_j)_{j\in\mathbb N}$ be an increasing sequence in $\mathbb N$
with nested spectra and with $\sup_{j\in\mathbb N}V(\nu_j)=\infty$. Then there
exists a function $f\in L_1[0,1)$ such that
\[
 \limsup_{j\to\infty}|s_{\nu_j}(f)(x)|=\infty
 \quad\text{for every }x\in[0,1),
\]
and a nonnegative function with the same property almost everywhere.
\end{corollary}

Without the nonnegative function, \cref{c7.5} is due to Goginava
and Oniani \cite[Theorem~2.1]{GO}.
\begin{proof}
Put $F\coloneqq \{\nu_j:j\in\mathbb N\}$. Let $K\in\mathbb N$ and choose $j_0$
with $V(\nu_{j_0})>K$. Let $j>j_0$. By \eqref{e7.4}, applied with $j_0$
and $j$ in place of $j$ and $j'$, the difference $\nu_j-\nu_{j_0}$ is
divisible by $2^{1+\max\operatorname{Sp}(\nu_{j_0})}$, and it is positive.
Moreover, $1\leq\nu_{j_0}<2^{1+\max\operatorname{Sp}(\nu_{j_0})}$ and
$V(\nu_{j_0})>K$. By \cref{d3.1}, applied with the quotient
$q\coloneqq (\nu_j-\nu_{j_0})/2^{1+\max\operatorname{Sp}(\nu_{j_0})}$ and the
remainder $a\coloneqq \nu_{j_0}$,
\[
 \frac{\nu_j-\nu_{j_0}}{2^{1+\max\operatorname{Sp}(\nu_{j_0})}}
 \in\mathcal H_{1+\max\operatorname{Sp}(\nu_{j_0}),K}(F)\qquad(j>j_0).
\]
These quotients increase strictly with $j$, and so
$\mathcal W_K(F)\geq\#\mathcal H_{1+\max\operatorname{Sp}(\nu_{j_0}),K}(F)=\infty$.
Since $K$ was arbitrary, \cref{t3.6,c3.5} give the two functions.
\end{proof}

\subsection{A convergent example}
If some of the lowest occupied positions of an integer are removed,
then its quotient does not change at the scales where a nonzero remainder
is left. The set in the following example is a union of finite sets, which
occupy disjoint intervals of positions, and the elements of each of these
sets arise in this way from its largest element. Along this union the
partial sums converge almost everywhere, but the union cannot be covered by
finitely many sets to which \cref{c7.2} or \cref{c7.3} applies.
\begin{example}\label{ex7.6}
For $j\in\mathbb N$ and $l\in\mathbb Z_+$ with $l<j$, define
\begin{align*}
 m_{j,l}\coloneqq \sum_{i=l}^{j-1}2^{2j^2+2i},\qquad
 F_j\coloneqq \{m_{j,l}:\ 0\leq l<j\},\qquad
 F\coloneqq \bigcup_{j=1}^{\infty}F_j.
\end{align*}
Then
\[
 \mathcal W_1(F)\leq1,\qquad \|s_*^F\|_{(1,1)}\lesssim1,\qquad
 \sup_{n\in F}V(n)=\infty.
\]
The partial sums along the increasing enumeration of $F$ converge almost
everywhere to $f$ for every $f\in L_1[0,1)$. Nevertheless, $F$ cannot be
covered by finitely many sets each having bounded variation or separated
spectra.
\end{example}
\begin{proof}
For any $j\in\mathbb N$ and $0\leq l<j$, we have
\[
 \operatorname{Sp}(m_{j,l})=\{2j^2+2i:\ l\leq i<j\}
 \subseteq[2j^2,2j^2+2j-2].
\]
Since $2j^2+2j-2<2(j+1)^2$, the intervals $[2j^2,2j^2+2j-2]$ are
pairwise disjoint. Each $F_j$ has $j$ distinct elements, so $F$ is
infinite. The $j-l$ ones in the binary expansion of $m_{j,l}$ are
separated by zeros. Hence
\[
 V(m_{j,l})=2(j-l),\qquad V(m_{j,0})=2j\longrightarrow\infty.
\]

Let $k\in\mathbb Z_+$. By \eqref{e7.3}, the quotient
$\lfloor m_{j,l}/2^k\rfloor$ and the remainder $m_{j,l}\bmod2^k$
are both positive if and only if
\[
 \min\operatorname{Sp}(m_{j,l})=2j^2+2l<k
 \leq2j^2+2j-2=\max\operatorname{Sp}(m_{j,l}).
\]
By the disjointness of the intervals $[2j^2,2j^2+2j-2]$, at most one
$F_j$ contains an element with both a positive quotient and a positive
remainder at scale $k$. For this $j$, the difference $m_{j,0}-m_{j,l}$ is a
sum of powers of two at positions strictly below $2j^2+2l<k$. Hence all
these elements have the same digits at the positions at least $k$, and
\[
 \left\lfloor m_{j,l}/2^k\right\rfloor
   =\left\lfloor m_{j,0}/2^k\right\rfloor.
\]
Thus $\#\mathcal H_{k,1}(F)\leq1$ for every $k$. \Cref{t3.4,t3.2} give the
maximal bound and the convergence; with $K=1$ and $N=1$, \eqref{e3.4}
gives $\|s_*^F\|_{(1,1)}\lesssim1$.

Finally, suppose that $F$ is covered by finitely many sets, each
having bounded variation or separated spectra. There are finitely many
sets with bounded variation, and so there is $K\in\mathbb Z_+$ such that
$V(n)\leq K$ for every element $n$ of these sets. Let $m\geq0$ be the
number of the sets with separated spectra. An element $n\in F$ with
$V(n)>K$ belongs to none of the sets with bounded variation, and hence to
one of the $m$ sets with separated spectra. By \eqref{e7.1}, the intervals
$[\min\operatorname{Sp}(n),\max\operatorname{Sp}(n)]$, where $n$ runs over a
set with separated spectra, are pairwise disjoint. All elements of $F_j$
have the highest occupied position $2j^2+2j-2$, so that each set with
separated spectra contains at most one element of $F_j$. On the other
hand, for $0\leq l<j$,
\[
 V(m_{j,l})=2(j-l)>K\iff l\leq j-\lfloor K/2\rfloor-1,
\]
and therefore
\[
 \#\{n\in F_j:\ V(n)>K\}=j-\lfloor K/2\rfloor\qquad(j>\lfloor K/2\rfloor).
\]
These elements belong to the $m$ sets with separated spectra, at most one
to each set. Taking $j=m+\lfloor K/2\rfloor+1$, we obtain
\[
 m+1=\#\{n\in F_j:\ V(n)>K\}\leq m,
\]
a contradiction.
\end{proof}

\subsection{A divergent example}
The increasing enumeration $(\nu_j)_{j\in\mathbb N}$ in \cref{ex7.6}
satisfies $\min\operatorname{Sp}(\nu_j)\to\infty$. The same limit holds
for every sequence with separated spectra. In the following example, we
construct a sequence with this limit and with
$\mathcal W_K(\{\nu_j:j\in\mathbb N\})=\infty$ for every $K\in\mathbb N$.

For any $f\in L_1[0,1)$ and $m,n\in\mathbb N$, the coefficient bound
$|\widehat f(r)|\leq\|f\|_{L_1[0,1)}$ gives
\[
 |s_n(f)(x)-s_m(f)(x)|\leq |n-m|\|f\|_{L_1[0,1)}
 \qquad(x\in[0,1)).
\]
Thus, as observed by Goginava and Oniani \cite[Remark~2.2]{GO},
unboundedness of the partial sums at a point passes to any sequence
that contains the given sequence as a subsequence or is at bounded
distance from it. Consequently, both conclusions of \cref{c7.5} hold
for every sequence with a subsequence at bounded distance from a
sequence with nested spectra and unbounded variation. The sequence of
the next example has no
subsequence at bounded distance from a sequence with nested spectra, so
that its divergence does not follow from \cref{c7.5} in this way.
\begin{example}\label{ex7.7}
The sequence $(\nu_j)_{j\in\mathbb N}$ defined by
\[
 \nu_j\coloneqq \sum_{i=j}^{2j}4^i\qquad(j\in\mathbb N)
\]
satisfies $\min\operatorname{Sp}(\nu_j)=2j\to\infty$, but there exists a
function $f\in L_1[0,1)$ such that
\[
 \limsup_{j\to\infty}|s_{\nu_j}(f)(x)|=\infty
 \quad\text{for every }x\in[0,1),
\]
and a nonnegative function with the same property almost everywhere.
There do not exist an increasing sequence $(m_l)_{l\in\mathbb N}$ in
$\mathbb N$ with nested spectra and a subsequence $(\nu_{j_l})_{l\in\mathbb N}$
such that
\[
 \sup_{l\in\mathbb N}|m_l-\nu_{j_l}|<\infty.
\]
In particular, no infinite subsequence of $(\nu_j)_{j\in\mathbb N}$ has nested
spectra.
\end{example}
\begin{proof}
Let $F\coloneqq \{\nu_j:j\in\mathbb N\}$. The sequence is strictly increasing, since
\[
 \nu_{j+1}-\nu_j=4^{2j+1}+4^{2j+2}-4^j>0.
\]
Fix $N\in\mathbb N$ and take $N\leq j\leq\lfloor3N/2\rfloor$. Division
by $2^{4N}=4^{2N}$ gives the remainder and the quotient
\begin{align*}
 \nu_j\bmod4^{2N}=\sum_{i=j}^{2N-1}4^i,\qquad
 \left\lfloor \nu_j/4^{2N}\right\rfloor
        =\sum_{i=0}^{2j-2N}4^i.
\end{align*}
The remainder has $2N-j$ isolated binary ones, so that
\[
 V\bigl(\nu_j\bmod4^{2N}\bigr)=2(2N-j)\geq N.
\]
The quotients are positive, and they are distinct, because their highest
occupied binary positions are $4j-4N$. For a fixed $K\in\mathbb N$
choose $N>K$. Then all these quotients belong to $\mathcal H_{4N,K}(F)$,
and
\[
 \#\mathcal H_{4N,K}(F)\geq\lfloor3N/2\rfloor-N+1
       =\lfloor N/2\rfloor+1\longrightarrow\infty\qquad(N\to\infty).
\]
Thus $\mathcal W_K(F)=\infty$ for every $K$, and \cref{t3.6,c3.5} give the
asserted functions.

Finally, suppose that a strictly increasing sequence $(m_l)_{l\in\mathbb N}$
in $\mathbb N$ has nested spectra and $|m_l-\nu_{j_l}|\leq C$ for all $l$, where $C\in\mathbb Z_+$ is fixed. The differences take only
finitely many integer values. After passing to a subsequence there is
a fixed integer $\delta$ with
\[
 m_l=\nu_{j_l}+\delta\qquad(l\in\mathbb N).
\]
By \cref{l7.4}, the new sequence $(m_l)_{l\in\mathbb N}$ still has
nested spectra. Choose an index $l_0$ so
large that $2^{1+\max\operatorname{Sp}(m_{l_0})}>|\delta|$; this is possible
because $m_l\to\infty$. By \eqref{e7.4},
\begin{align}\label{e7.6}
 m_l\bmod2^{1+\max\operatorname{Sp}(m_{l_0})}=m_{l_0}\qquad(l\geq l_0).
\end{align}
For all sufficiently large $l$, the integer $\nu_{j_l}$ is divisible by
$2^{1+\max\operatorname{Sp}(m_{l_0})}$, since its lowest occupied position
is $2j_l\to\infty$. Fix such an $l\geq l_0$. Since
$m_l=\nu_{j_l}+\delta$, we have
\[
 m_{l_0}-\delta=\nu_{j_l}-(m_l-m_{l_0}),
\]
and both terms on the right-hand side are divisible by
$2^{1+\max\operatorname{Sp}(m_{l_0})}$, the first one by the choice of $l$
and the second one by \eqref{e7.6}. Moreover, $m_{l_0}-\delta=\nu_{j_{l_0}}$,
$m_{l_0}<2^{1+\max\operatorname{Sp}(m_{l_0})}$ and
$|\delta|<2^{1+\max\operatorname{Sp}(m_{l_0})}$, so that
\[
 0<\nu_{j_{l_0}}=m_{l_0}-\delta<2\cdot2^{1+\max\operatorname{Sp}(m_{l_0})}.
\]
Thus $\nu_{j_{l_0}}$ is a positive multiple of
$2^{1+\max\operatorname{Sp}(m_{l_0})}$ smaller than twice this power, and so
$\nu_{j_{l_0}}=2^{1+\max\operatorname{Sp}(m_{l_0})}$. This is impossible,
since every $\nu_j$ has $j+1\geq2$ occupied binary positions. The
assertion follows.
\end{proof}

\section*{Declarations}

\noindent {\bf Ethical approval} Not applicable. \\

\noindent {\bf Competing interests} The authors declare that there are no competing interests. \\

\noindent {\bf Funding} This project is partially supported by the National Natural
Science Foundation of China (Grant Nos. 12671164, W2411005, 12431006, and 12371093)
and the Natural Science Foundation of Hunan Province of China (Grant No.
2025JJ40007). \\

\noindent {\bf Availability of data and materials} Data sharing not applicable to
this article as no datasets were generated or analysed during the current study.\\


\begin{thebibliography}{99}

\bibitem{APT25}
N.~Areshidze, L.-E.~Persson and G.~Tephnadze,
\emph{Convergence almost everywhere of partial sums and F\'ejer means
of Vilenkin--Fourier series},
Mediterr. J. Math. \textbf{22} (2025), Paper No.~15, 17 pp.

\bibitem{Billard}
P.~Billard,
\emph{Sur la convergence presque partout des s\'eries de Fourier--Walsh
des fonctions de l'espace $L^2(0,1)$},
Studia Math. \textbf{28} (1967), 363--388.

\bibitem{BG26}
I.~Blahota and G.~G\'at,
\emph{Almost everywhere divergence of Ces\`aro means of subsequences of
Walsh--Paley Fourier partial sums},
preprint (2026),
arXiv:2606.13879.

\bibitem{Carleson}
L.~Carleson,
\emph{On convergence and growth of partial sums of Fourier series},
Acta Math. \textbf{116} (1966), 135--157.

\bibitem{Demeter}
C.~Demeter,
\emph{A guide to Carleson's theorem},
Rocky Mountain J. Math. \textbf{45} (2015), 169--212.

\bibitem{DP}
F.~Di Plinio,
\emph{Lacunary Fourier and Walsh--Fourier series near $L^1$},
Collect. Math. \textbf{65} (2014), 219--232.

\bibitem{DL}
Y.~Q.~Do and M.~T.~Lacey,
\emph{On the convergence of lacunary Walsh--Fourier series},
Bull. Lond. Math. Soc. \textbf{44} (2012), 241--254.

\bibitem{Fine}
N.~J.~Fine,
\emph{On the Walsh functions},
Trans. Amer. Math. Soc. \textbf{65} (1949), 372--414.

\bibitem{Gat10}
G.~G\'at,
\emph{Almost everywhere convergence of Fej\'er and logarithmic means of
subsequences of partial sums of the Walsh--Fourier series of integrable
functions},
J. Approx. Theory \textbf{162} (2010), 687--708.

\bibitem{Gat19}
G.~G\'at,
\emph{On the convergence of Fej\'er means of some subsequences of partial
sums of Walsh--Fourier series},
Ann. Univ. Sci. Budapest. Sect. Comput. \textbf{49} (2019), 187--198.

\bibitem{Gog26}
U.~Goginava,
\emph{Almost everywhere convergence of arithmetic means of Walsh--Fourier
partial sums along subsequences},
preprint (2026),
arXiv:2605.05146.

\bibitem{GM26}
U.~Goginava and F.~Mukhamedov,
\emph{Dyadic martingale transforms and weighted Walsh--Carleson operators},
preprint (2026),
arXiv:2605.07667.

\bibitem{GO}
U.~Goginava and G.~G.~Oniani,
\emph{On the divergence of subsequences of partial Walsh--Fourier sums},
J. Math. Anal. Appl. \textbf{497} (2021), Paper No.~124900.

\bibitem{Gosselin}
R.~P.~Gosselin,
\emph{On the divergence of Fourier series},
Proc. Amer. Math. Soc. \textbf{9} (1958), 278--282.

\bibitem{HW78}
D.~C.~Harris and W.~R.~Wade,
\emph{Sets of divergence on the group $2^\omega$},
Trans. Amer. Math. Soc. \textbf{240} (1978), 385--392.

\bibitem{Hunt}
R.~A.~Hunt,
\emph{On the convergence of Fourier series},
in \emph{Orthogonal Expansions and their Continuous Analogues}
(Proc. Conf., Edwardsville, Ill., 1967),
Southern Illinois University Press, Carbondale, Ill., 1968, pp.~235--255.

\bibitem{Katz66}
Y.~Katznelson,
\emph{Sur les ensembles de divergence des s\'eries trigonom\'etriques},
Studia Math. \textbf{26} (1966), 301--304.

\bibitem{K93}
S.~V.~Konyagin,
\emph{The Fourier--Walsh subsequence of partial sums},
Math. Notes \textbf{54} (1993), 1026--1030.

\bibitem{K05}
S.~V.~Konyagin,
\emph{Divergence everywhere of subsequences of partial sums of
trigonometric Fourier series},
Proc. Steklov Inst. Math. (2005), suppl.~2, S167--S175.

\bibitem{K06}
S.~V.~Konyagin,
\emph{Almost everywhere convergence and divergence of Fourier series},
in \emph{International Congress of Mathematicians, Madrid 2006},
Vol.~II, European Mathematical Society, Z\"urich, 2006, pp.~1393--1403.

\bibitem{LaceyThiele}
M.~T.~Lacey and C.~Thiele,
\emph{A proof of boundedness of the Carleson operator},
Math. Res. Lett. \textbf{7} (2000), 361--370.

\bibitem{Lie12}
V.~Lie,
\emph{On the pointwise convergence of the sequence of partial
Fourier sums along lacunary subsequences},
J. Funct. Anal. \textbf{263} (2012), 3391--3411.

\bibitem{Lie13}
V.~Lie,
\emph{On the boundedness of the Carleson operator near $L^1$},
Rev. Mat. Iberoam. \textbf{29} (2013), 1239--1262.

\bibitem{Lie}
V.~Lie,
\emph{Pointwise convergence of Fourier series (I). On a conjecture of Konyagin},
J. Eur. Math. Soc. \textbf{19} (2017), 1655--1728.

\bibitem{Lie19}
V.~Lie,
\emph{The pointwise convergence of Fourier series (II).
Strong $L^1$ case for the lacunary Carleson operator},
Adv. Math. \textbf{357} (2019), Paper No.~106831.

\bibitem{Luk95}
S.~F.~Lukomskij,
\emph{A criterion for the almost-everywhere convergence of
Fourier--Walsh square partial sums of integrable functions},
Sb. Math. \textbf{186} (1995), 1057--1070.

\bibitem{Luk98}
S.~F.~Lukomskij,
\emph{Multiple Walsh series: Convergence in measure and almost everywhere},
Dokl. Math. \textbf{57} (1998), 81--82.

\bibitem{NOT}
F.~Nazarov, R.~Oberlin and C.~Thiele,
\emph{A Calder\'on Zygmund decomposition for multiple frequencies
and an application to an extension of a lemma of Bourgain},
Math. Res. Lett. \textbf{17} (2010), 529--545.

\bibitem{OT}
R.~Oberlin and C.~Thiele,
\emph{New uniform bounds for a Walsh model of the bilinear Hilbert transform},
Indiana Univ. Math. J. \textbf{60} (2011), 1693--1712.

\bibitem{On22}
G.~G.~Oniani,
\emph{On the divergence sets of Fourier series in systems of characters
of compact Abelian groups},
Math. Notes \textbf{112} (2022), 100--108.

\bibitem{On23}
G.~G.~Oniani,
\emph{Almost everywhere convergence of nets of operators and weak
type maximal inequalities},
Fund. Math. \textbf{262} (2023), 287--304.

\bibitem{Paley}
R.~E.~A.~C.~Paley,
\emph{A remarkable series of orthogonal functions (I)},
Proc. London Math. Soc. (2) \textbf{34} (1932), 241--264.

\bibitem{wk6}
L.-E.~Persson, G.~Tephnadze and F.~Weisz,
\emph{Martingale Hardy Spaces and Summability of One-Dimensional
Vilenkin--Fourier Series},
Birkh\"auser, Cham, 2022.

\bibitem{Sawyer}
S.~Sawyer,
\emph{Maximal inequalities of weak type},
Ann. of Math. (2) \textbf{84} (1966), 157--174.

\bibitem{Schipp67}
F.~Schipp,
\emph{\"Uber die Gr\"o\ss enordnung der Partialsummen der Entwicklung
integrierbarer Funktionen nach W-Systemen},
Acta Sci. Math. (Szeged) \textbf{28} (1967), 123--134.

\bibitem{Schipp69}
F.~Schipp,
\emph{\"Uber die Divergenz der Walsh--Fourierreihen},
Ann. Univ. Sci. Budapest. E\"otv\"os Sect. Math. \textbf{12} (1969), 49--62.

\bibitem{SWS}
F.~Schipp, W.~R.~Wade and P.~Simon, with the assistance of J.~P\'al,
\emph{Walsh Series: An Introduction to Dyadic Harmonic Analysis},
Adam Hilger, Bristol and New York, 1990.

\bibitem{Sjolin}
P.~Sj\"olin,
\emph{An inequality of Paley and convergence a.e. of Walsh--Fourier series},
Ark. Mat. \textbf{7} (1969), 551--570.

\bibitem{Stein}
E.~M.~Stein,
\emph{On limits of sequences of operators},
Ann. of Math. (2) \textbf{74} (1961), 140--170.

\bibitem{Thiele00}
C.~Thiele,
\emph{The quartile operator and pointwise convergence of Walsh series},
Trans. Amer. Math. Soc. \textbf{352} (2000), 5745--5766.

\bibitem{Totik}
V.~Totik,
\emph{On the divergence of Fourier series},
Publ. Math. Debrecen \textbf{29} (1982), 251--264.

\bibitem{Walsh}
J.~L.~Walsh,
\emph{A closed set of normal orthogonal functions},
Amer. J. Math. \textbf{45} (1923), 5--24.

\bibitem{wk}
F.~Weisz,
\emph{Martingale Hardy Spaces and their Applications in Fourier
Analysis},
Lecture Notes in Math. \textbf{1568}, Springer, Berlin, 1994.

\bibitem{wk2}
F.~Weisz,
\emph{Summability of Multi-Dimensional Fourier Series and Hardy Spaces},
Mathematics and Its Applications \textbf{541}, Kluwer Academic Publishers,
Dordrecht, Boston, London, 2002.

\end{thebibliography}
\end{document}